\documentclass[12pt]{article}

\usepackage{amssymb,amsmath,amsthm, amsfonts}

\newcommand{\To}{\longrightarrow}

\newcommand{\f}{\mathcal{F}}
\newcommand{\e}{\mathcal{E}}

\def\sqr#1#2{{\vcenter{\vbox{\hrule height.#2pt
              \hbox{\vrule width.#2pt height#1pt \kern#1pt \vrule width.#2pt}
              \hrule height.#2pt}}}}
\def\signed #1{{\unskip\nobreak\hfil\penalty50
              \hskip2em\hbox{}\nobreak\hfil#1
              \parfillskip=0pt \finalhyphendemerits=0 \par}}
\def\endpf{\signed {$\sqr69$}}
\def\3n{\negthinspace \negthinspace \negthinspace }
\def\2n{\negthinspace \negthinspace }
\def\1n{\negthinspace }

\def\={\buildrel \triangle \over =}

\def\a{\alpha}
\def\b{\beta}

\def\d{\delta}
\def\e{\varepsilon}

\def\l{\lambda}

\def\si{\sigma}

\def\th{\theta}
\def\o{\omega}

\def\u{\upsilon}

\def\L{\Lambda}

\def\F{\Phi}

\def\cF{{\cal F}}

\def\cL{{\cal L}}

\def\cP{{\cal P}}

\def\ms{\medskip}

\def\q{\quad}

\def\limsup{\mathop{\overline{\rm lim}}}
\def\liminf{\mathop{\underline{\rm lim}}}

\def\esssup{\mathop{\rm esssup}}
\def\max{\mathop{\rm max}}
\def\min{\mathop{\rm min}}
\def\exp{\mathop{\rm exp}}
\def\sup{\mathop{\rm sup}}

\def\inf{\hbox{\rm inf$\,$}}
\def\esssup{\hbox{\rm ess$\,$\rm sup$\,$}}

\def\|{\Big |}
\def\({\Big (}
\def\){\Big )}
\def\[{\Big[}
\def\]{\Big]}
\def\be{\begin{equation}}
\def\bel{\begin{equation}\label}
\def\ee{\end{equation}}
\def\bea{\begin{eqnarray}}
\def\eea{\end{eqnarray}}
\def\bea*{\begin{eqnarray*}}
\def\eea*{\end{eqnarray*}}
\def\bt{\begin{theorem}}
\def\bcd{\begin{condition}}
\def\ecd{\end{condition}}
\def\et{\end{theorem}}
\def\bc{\begin{corollary}}
\def\ec{\end{corollary}}
\def\bde{\begin{definition}}
\def\ede{\end{definition}}
\def\bl{\begin{lemma}}
\def\el{\end{lemma}}
\def\bp{\begin{proposition}}
\def\ep{\end{proposition}}
\def\br{\begin{remark}}
\def\er{\end{remark}}
\def\ba{\begin{array}}
\def\ea{\end{array}}
\def\ed{\end{document}}

\def\square#1{\vbox{\hrule\hbox{\vrule height#1%
     \kern#1\vrule}\hrule}}
\def\rectangle#1#2{\vbox{\hrule\hbox{\vrule height#1%
     \kern#2\vrule}\hrule}}

\def\essinf{{\rm { ess.inf}}}
\def\esssup{{\rm ess.sup}}

\font\tenbb=msbm10 \font\sevenbb=msbm7 \font\fivebb=msbm5

\newfam\bbfam
\scriptscriptfont\bbfam=\fivebb \textfont\bbfam=\tenbb
\scriptfont\bbfam=\sevenbb

\newtheorem{lemma}{Lemma}[section]
\newtheorem{remark}{Remark}[section]
\newtheorem{example}{Example}[section]
\newtheorem{theorem}{Theorem}[section]
\newtheorem{corollary}{Corollary}[section]

\newtheorem{definition}{Definition}[section]
\newtheorem{proposition}{Proposition}[section]
\newtheorem{condition}{Condition}[section]

\makeatletter
   
   \@addtoreset{equation}{section}
\makeatother

\begin{document}

    \title{{\bf  Backward SDEs with superquadratic growth}}

\author{Freddy Delbaen\thanks{ Department of Mathematics, ETH Z\"{u}rich,
Switzerland. {\small\it E-mail:} {\small\tt
 delbaen@math.ethz.ch}.
 Part of the research was done while this author
was visiting Princeton, Vancouver and China. The hospitalities of 
Princeton University, UBC, Fudan University and Shandong University are greatly
appreciated. \ms},
\quad Ying Hu\thanks{ IRMAR, Universit\'e Rennes 1, 35 042 RENNES
Cedex, France. {\small\it E-mail:} {\small\tt
 ying.hu@univ-rennes1.fr}.\ms} \quad and\quad Xiaobo Bao\thanks{Department of Mathematics, ETH Z\"{u}rich,
 Switzerland. {\small\it E-mail:} {\small\tt
 baoxb@math.ethz.ch}.\ms }
}

\date{February 19, 2009}
\maketitle

\abstract{In this paper, we discuss the solvability of backward
stochastic differential equations (BSDEs) with superquadratic
generators. We first prove  that given a superquadratic generator,
there exists a bounded terminal value, such that the associated BSDE
does not admit any bounded solution. On the other hand, we prove
that if the superquadratic BSDE admits a bounded solution, then
there exist infinitely many bounded solutions for this BSDE.
Finally, we prove the existence of a solution for Markovian BSDEs
where the terminal value is a bounded continuous function of a
forward stochastic differential equation.
 }

\section{Introduction.}
Since the pioneer works on BSDEs of Bismut \cite{Bis} and Pardoux-Peng
 \cite{PP90}, lots of works
have been done in this area and the original Lipschitz assumption on
the generator, i.e., the function $g$ in the BSDE:
\begin{equation}\label{bsde}
Y_t=\xi-\int_t^T g(s,Y_s,Z_s)\,ds+\int_t^T Z_s\,dB_s,\quad 0\le t\le
T,
\end{equation}
has been weakened in many situations. Let us recall that, in the
previous BSDE, we are looking for a pair of processes $(Y,Z)$ which
is required to be  predictable with respect
to the filtration generated by the Brownian motion $B$. One of the
most important works in this direction is that of Kobylanski
\cite{Kobylanski} concerning scalar-valued quadratic BSDEs with
bounded terminal value. We should point out that quadratic BSDE
means a BSDE whose generator has at most a quadratic growth with
respect to the variable $z$. For these quadratic BSDEs, all the
classical results, existence and uniqueness, comparison and
stability of solutions, have been stated in \cite{Kobylanski} but
with the restriction that the terminal conditions have to be bounded
random variables. Recently, existence and uniqueness of solutions of
quadratic BSDEs with unbounded terminal value were studied by Briand
and Hu in \cite{BH1,BH2}.

In this paper, we study the solvability of superquadratic BSDE
(\ref{bsde}) whose generator $g$ is superquadratic, i.e.,
$$\limsup_{|z|\rightarrow +\infty} \frac{g(z)}{|z|^2}=\infty.$$
We shall study this BSDE with bounded terminal value. And in
addition, we suppose that $g$ is a deterministic convex (or concave)
function which is independent of $y$ with $g(0)=0$.

The first part of this paper  shows the ill-posedness of these
BSDEs. We first prove that given a superquadratic generator, there
always exists a bounded terminal value, such that the associated
BSDE does not admit any bounded solution. On the other hand, we
prove that if the superquadratic BSDE admits a bounded solution,
then there exist infinitely many bounded solutions for this BSDE.
And finally, we show that the monotone stability, which plays a
crucial role in quadratic BSDEs (see, e.g., \cite{Kobylanski,BH1}),
does not hold.

\bigskip

In the second part of this paper, we study BSDE (\ref{bsde}) in the
Markovian case, i.e., the terminal value
$$\xi=\Phi(X_T^{t,x}),$$
where the diffusion process $X$ is the solution to the SDE:
\begin{equation}\label{sde}
X_s=x+\int_t^s b(r,X_r)\,dr+\int_t^s \sigma \,dB_r,\quad t\le s\le
T.
\end{equation}

It is by now well-known (see, e.g., \cite{PP92,Kobylanski,BH2} ) that, if $g$ is Lipschitz or quadratic, there exists a
link between the solution of (\ref{bsde}) and that of the following
PDE: \be~\label{pde} \left\{\ba{l}
u_t(t,x)+\frac{1}{2}\mbox{trace}\big(\sigma\sigma^Tu_{xx}(t,x)\big)+u_x(t,x)b(t,x)
-g(-u_x(t,x)\sigma)=0,\\u(T,x)=\Phi(x).
    \ea\right.
\ee

 This type of PDE (called viscous Hamilton-Jacobi equation)
is already well studied when $\sigma$ is the identity and
$g(z)=-|z|^p$, see, e.g., Gilding et al. \cite{ggk} and Ben-Artzi et al.
\cite{bsw}. In particular, in \cite{ggk}, they established the
existence and uniqueness of classical solution to this PDE when
$\sigma$ is the identity.

We prove that in the Markovian case, the BSDE (\ref{bsde}) admits a
solution when $\Phi$ is bounded and continuous. Moreover, if we
define
$$u(t,x)=Y_t^{t,x},$$
then $u$ is a continuous viscosity solution to PDE (\ref{pde}). We
note that in our case, some kind of degeneracy of $\sigma$ is
allowed, whereas in \cite{ggk} and \cite{bsw}, they assumed that
$\sigma$ is the identity.

A key idea to prove the existence in the Markovian case comes from
the following a priori estimate of $Z$:
$$|Z_t|\le c ||\Phi||_\infty (T-t)^{-\frac{1}{2}},$$
where $c>0$ is a constant. We prove this inequality by using a
stochastic argument based on BMO martingales and Jensen's
inequality. Note that Gilding et al. \cite{ggk} proved the same type
of a priori estimate for $u_x$ when $\sigma$ is identity, by use of
Bernstein's method.

The paper is organized as follows: in the next section, we give some
preliminaries about the connection between dynamic utility functions
and BSDEs. Section 3 shows the ill-posedness in the general case.
The last section is devoted to the proof of the existence of a
solution in the Markovian case.

\section{Dynamic Utility Functions and Backward SDEs. }
Let $\{B_t,0 \leq t \leq T\}$ be a $d$-dimensional standard Brownian
motion defined on a probability space $(\Omega,\mathcal{F},P)$. Let
$\{\mathcal{F}_t,0 \leq t \leq T\}$ be the natural filtration of
$\{B_t, t\in [0,T]\}$, augmented by all $P$-null sets of
$\mathcal{F}$.

 Before recalling the definition of dynamic utility
functions, we need the following notations.
$$\ba{rcl} L^{\infty}(\mathcal{F}_T)&:=&\{\xi\,: \hbox{ \rm  bounded and $\mathcal{F}_T$-measurable random variable }\},\\
 \cL_{\cF}^2(0,T;{\mathbb{R}}^m)&:=&\{\varphi: \hbox{ \rm $\mathbb{R}^m$-valued,  $\{\cF_t\}_{ 0\le t\le
T}$-predictable
 and }  E\left[\int_0^T|\varphi_t|^2\,dt\right]<\infty\}.\ea$$

We identify random variables that are equal $P$ a.s.
 \bde We call a dynamic utility function with the Fatou property any family of
operators, indexed by stopping times $\sigma$
$$U_{\sigma}:L^{\infty}(\cF_T)\rightarrow L^{\infty}(\cF_{\sigma})$$
and satisfying:
\begin{itemize}
    \item (A1) {\bf Positivity:}  $U_{\sigma}(0)=0, U_{\sigma}(\xi)\geq 0 \mbox{ for
    all } \xi \geq 0.$
    \item (A2) {\bf Concavity:} $U_{\sigma}(t\xi +(1-t)\eta)\geq tU_{\sigma}(\xi)+(1-t)U_{\sigma}(\eta), \mbox{ for all }t,
    0\leq t \leq 1 \mbox{ and all }\xi, \eta \in L^{\infty}.$
    \item (A3) {\bf Translability:} $U_{\sigma}(\xi +a)=U_{\sigma}(\xi)+ a, \mbox{ for all }a \in L^{\infty}(\cF_{\sigma}).$
    \item (A4) {\bf Fatou property:} Given a sequence $(\xi_n)_{n\geq
    1}$, such that $\sup||\xi_n||_{\infty} < \infty$, then
    $\xi_n\downarrow\xi$ a.s. implies $U_{\sigma}(\xi)=
    {\lim}_{n\rightarrow \infty}U_{\sigma}({\xi_n})$ a.s.
\end{itemize}    \ede
For a lower semi-continuous convex function $f:R^d\rightarrow
R_{+}\cup\{\infty\}$ such that $f(0)=0$ and for $\xi \in
L^{\infty}(\cF_T)$, we define \be~\label{ut}
U_{\sigma}(\xi)=\essinf\bigg\{E_{Q}\Big[\xi+\int_{\sigma}^Tf(q_u)\,du\Big|\f_{\sigma}\Big]\,\bigg|\,
Q\sim P\bigg\},\ee where $\si\in [0,T]$ is a stopping time and the
density process $E_{P}[\frac{dQ}{dP}|\f_t]=\mathcal{E}(q\cdot
B)_t=\exp(\int_0^t q_u \,dB_u-\frac{1}{2}\int_0^t|q_u|^2\,du)$. It
is easy to prove that $U$ is a dynamic utility function. As shown by
Delbaen-Peng-Rosazza Gianin \cite{DPR}, $U$ is time consistent and
all time consistent dynamic utility functions are of a similar form.

Set $C_0(Q)=E_Q\big[\int_0^Tf(q_u)\,du\big]$ and $\cP=\{Q\mid Q \ll
P \}.$ The utility function $U_0$ can be defined by $\cP$.

\bl  ~\label{qplemma}For any $\xi \in L^{\infty}(\cF_T)$, \be
~\label{qp}
U_{0}(\xi)=\inf\bigg\{E_{Q}\Big[\xi+\int_{0}^Tf(q_u)\,du\Big]\,\bigg|\,
Q\in \cP\bigg\}.\ee \el

\emph{Proof.} For any  $Q\in \cP$ with
$L_t=E_P\big[\frac{dQ}{dP}\big|\f_t\big]=\mathcal{E}(q\cdot B)_t$,
using It\^o's lemma we get that the density process of
$Q_{\lambda}\triangleq \lambda Q+ (1-\l)P$ is
$\mathcal{E}(q_{\l}\cdot B)$ with
$$q_{\l}(t)=\frac{\l L_t q_t}{\l L_t+(1-\l)}1_{\{t\leq \tau\}},$$
where $\tau=\inf\{t\in [0,T]\ |\ L_t=0\}\wedge T$ is a stopping
time.\\ Then from the convexity of $f$: \bea* C_0(Q_{\l})&=& E_{Q_{\l}}\Big[\int_0^Tf(q_{\l}(u))\,du\Big]\\
&\leq& E_{Q_{\l}}\Big[\int_0^{\tau}\frac{\l L_t}{\l
L_t+(1-\l)}f(q(t))\,dt\Big] \\
&=& E_P\Big[\int_0^{\tau}{\l L_t}f(q(t))\,dt\Big]\\
&=& E_Q\Big[\int_0^{\tau}{\l }f(q(t))\,dt\Big] \\
&=&\l C_0(Q),\eea* where $\l \in [0,1]$, we deduce that
$\overline{\lim}_{\l \rightarrow 1}C_0(Q_{\l})\leq C_0(Q).$ \\
Notice that for any $\l \in [0,1)$, $Q_{\l}$ is equivalent to $P$.
Thus
$$\inf\bigg\{E_{Q}\Big[\xi+\int_{0}^Tf(q_u)\,du\Big]\,\bigg|\, Q\in
\cP\bigg\}\geq
\inf\bigg\{E_{Q}\Big[\xi+\int_{0}^Tf(q_u)\,du\Big]\,\bigg|\,
Q\sim P\bigg\}.$$ Since $\{Q\,|\,Q\sim P\} \subseteq \{Q\mid Q\ll
P\}$, we have
$$U_0(\xi)=\inf\bigg\{E_{Q}\Big[\xi+\int_{0}^Tf(q_u)\,du\Big]\,\bigg|\, Q\in
\cP\bigg\}.$$ \endpf

\br The function $C_0:\cP\rightarrow \overline{R_+}$ is lower
semi-continuous (just use Fatou's lemma) and convex. A duality
argument then shows that for $Q\in \cP$
$$C_0(Q)=\sup\left\{E_Q[-\xi]\ \Big| \  U_0(\xi)\geq 0\right\}.$$
In other words $C_0$ is the minimal penalty function as defined in
F\"{o}llmer-Schied \cite{FS}. We also remark that for $Q\ll P$, the
previous reasoning and the lower semi-continuity imply
$C_0(Q_{\lambda})\rightarrow C_0(Q)$.

\er

 However, for a stopping time $\sigma$, $U_{\sigma}(\xi)$ cannot
be the essential infimum over $\cP$ $P$ a.s. Instead, by the similar
technique as that in Lemma~\ref{qplemma}, we have:

\br ~\label{utq} For any measure $Q^*\in\cP$ and $\xi \in
L^{\infty}(\cF_T)$, \be
U_{\sigma}(\xi)=\essinf\bigg\{E_{Q}\Big[\xi+\int_{\sigma}^Tf(q_u)\,du\Big|\f_{\sigma}\Big]\,\bigg|\,
Q\in \cP, Q\sim P \mbox{ on } \f_{\sigma}\bigg\},\ \ P\ a.s., \ee
for any stopping time $\sigma\in [0,T]$ and, \be
U_{\sigma}(\xi)=\essinf\bigg\{E_{Q}\Big[\xi+\int_{\sigma}^Tf(q_u)\,du\Big|\f_{\sigma}\Big]\,\bigg|\,
Q\in \cP, Q^*\ll Q\bigg\},\ \ Q^*\ a.s.\ee \er

 \bp ~\label{qmartin}For any $\xi\in L^{\infty}(\f_T)$, the dynamic utility function $U$
 defined by (\ref{ut}) has the following properties:\\
 1) For all $Q\ll P$, we have that $U_t(\xi)+\int_0^{\tau\wedge t} f(q_u)\,du$ is a
 $Q$-submartingale where $\tau=\inf\{t\in [0,T]\ |\ L_t=0\}$.\\
  2) If there is a probability
 measure
 $Q\ll P$ with $U_0(\xi)=E_{Q}[\xi+\int_0^{\tau}f(q_u)\,du]$, then $U_t(\xi)+\int_0^{\tau\wedge t } f(q_u)\,du$ is a
 $Q$-martingale.
\ep

 \emph {Proof.} 1) For any $s<t$, it follows from Remark
 \ref{utq} that for any $Q\ll P$,
  \bea*
 &&E_Q\left[U_t(\xi)+\int_{\tau\wedge s}^{\tau\wedge
t}f(q_u)\,du\bigg|\f_s\right]\\
&=&E_{Q}\left[\(\essinf_{Q'\sim P}
 \Big\{E_{Q^{'}}\big[\xi+\int_{t}^{T}f(q'_u)\,du\big|\f_t\big]\Big\}+
 \small\int_{\tau\wedge s}^{\tau\wedge t}f(q_u)\,du\)\bigg|\f_s\right] \\
 &\geq&\essinf\bigg\{E_{Q^{''}}\Big[\xi+\int_{\tau\wedge s}^{T} f(q''_u)\,du\Big|\f_s\Big] \,\bigg |\,q''_u=q'_u+1_{\{\tau\wedge s\leq u\leq t\}}(q_u-q'_u)\bigg\}\\
&\geq&\essinf\bigg\{E_{Q^{''}}\Big[\xi+\int_{\tau\wedge s}^{T} f(q''_u)\,du\Big|\f_s\Big] \,\bigg |\,Q''\in\cP,Q\ll Q''\bigg\}\\
 &\geq& U_s(\xi),\ \ Q\ \mbox{\,a.s.}
 \eea*
Hence, $$U_s(\xi)+\int_0^{\tau \wedge s} f(q_u)\,du\leq
E_Q\Big[U_t(\xi)+\int_0^{\tau \wedge t} f(q_u)\,du\Big|\f_s\Big],
\,Q\ a.s.$$ Therefore, we have $U_t(\xi)+\int_0^{\tau\wedge t}
f(q_u)\,du$ is a
 $Q$-submartingale. \\
2) As $Q$ is absolutely continuous with respect to $P$, it follows
from the result we just proved, that \be ~\label{subm} U_0(\xi) \leq
E_Q\Big[U_t(\xi)+\int_0^{\tau\wedge t} f(q_u)\,du\Big].\ee Combining
$U_0(\xi)=E_{Q}[\xi+\int_0^{\tau}f(q_u)\,du]$ with the inequality
(\ref{subm}), we have
$$E_Q\Big[\xi+\int_{\tau \wedge t}^{\tau}f(q_u)\,du\Big]\leq E_Q[U_t(\xi)].$$
This implies that \be U_t(\xi)=E_Q\Big[\xi+\int_{\tau \wedge
t}^{\tau}f(q_u)\,du\Big|\f_t\Big],\ Q\ \mbox{a.s.}\ee Thus
$U_t(\xi)+\int_0^{\tau\wedge t }
f(q_u)\,du=E_Q[\xi+\int_{0}^{\tau}f(q_u)\,du|\f_t]$ is a $Q$-
martingale.\endpf

\br In the above proposition, $\tau$ can be replaced by $T$ since
$Q[\tau=T]=1$.\er

 \br In particular, we have that the process $\{U_t(\xi),\
t\in [0,T]\}$ is a $P$-submartingale. Thus there exists a
c\`{a}dl\`{a}g version. \er

 For any $\xi \in L^{\infty}(\f_T)$,  $|U_t(\xi)|\leq \parallel\xi\parallel_{\infty}$. So applying the Doob-Meyer decomposition theorem, there exists a unique nondecreasing
predictable process $\{A_t\}_{0\leq t\leq T}$ with $A_0=0$ and a
 continuous martingale $\{M_t\}_{0\leq t\leq T}$ with $M_0=0$, such that
\be U_t(\xi)=U_0(\xi)+A_t-M_t.\ee
 \bl\label{bmo} For all
$\xi \in L^{\infty}(\f_T)$, the martingale part $\{M_t\}_{0\leq
t\leq T}$ of $U(\xi)$ induced by the Doob-Meyer decomposition
theorem is a BMO-martingale. \el

\emph {Proof.} For a given $\xi \in L^{\infty}(\f_T)$,
$|U_t(\xi)|\leq \parallel\xi\parallel_{\infty}$. Then applying
It\^o's formula to $(U_t(\xi)+\parallel\xi\parallel_{\infty})^2$, we
get\bea*
&&(U_t(\xi)+\parallel\xi\parallel_{\infty})^2+\int_t^Td\langle
M,M\rangle_s\\&=&(\xi+\parallel\xi\parallel_{\infty})^2-2\int_t^T(U_{s-}(\xi)+\parallel\xi\parallel_{\infty})\,dA_s
-\int_t^TdK_s\\
&&+2\int_t^T(U_{s-}(\xi)+\parallel\xi\parallel_{\infty})\,dM_s,
\eea* where \begin{eqnarray*}K_s&:=&\sum_{r\leq
s}\Big\{(U_r(\xi)+\parallel\xi\parallel_{\infty})^2-(U_{r-}(\xi)+\parallel\xi\parallel_{\infty})^2\\
&&-2(U_{r-}(\xi)+\parallel\xi\parallel_{\infty})(U_{r}(\xi)-U_{r-}(\xi))\Big\}\\
&=&\sum_{r\leq s}\(U_r(\xi)-U_{r-}(\xi)\)^2
\end{eqnarray*}
is an increasing process. Hence,
$$(U_t(\xi)+\parallel\xi\parallel_{\infty})^2+\int_t^Td\langle M,M\rangle_s\,\leq(\xi+\parallel\xi\parallel_{\infty})^2
+2\int_t^T(U_{s-}(\xi)+\parallel\xi\parallel_{\infty})\,dM_s$$ from
which we deduce, for any stopping time $0\leq\si\leq T$,
$$E\left[\int_{\si}^Td\langle M,M\rangle_t\|\f_{\si}\right]\leq 4\parallel\xi\parallel^2_{\infty}.$$

 Therefore, $\parallel
M\parallel_{BMO_2}\leq 2\parallel\xi\parallel_{\infty} $ which
completes the proof. \endpf

The predictable representation theorem implies that there exists a
predictable process $Z\in\cL_{\cF}^2(0,T;{\mathbb{R}}^d)$ such that
\be~\label{mr} M_t=\int_0^tZ_s\,dB_s.\ee So we get
\be~\label{decom-az} U_t(\xi)=U_0(\xi)+A_t-\int_0^tZ_s\,dB_s.\ee

If $g: R^{ d}\rightarrow R_{+}\cup \{\infty\}$ is the
Fenchel-Legendre transform of $f$:
$$g(z)=\sup_{x\in R^d}(zx-f(x)),$$
then $g$ is also convex and $g(0)=0.$

We make the standard assumption such that \emph{both $f$ and $g$ are
finite}. We do not treat the case where $f$ or $g$ can take the
value $+\infty$. This case is similar and only requires cosmetic
changes. To make the paper simpler, we dropped this more general
case.

 \bt ~\label{atgz} Let $U$ be the dynamic utility
function defined by (\ref{ut}) and let
$U_0(\xi)+A_t-\int_0^tZ_udB_u$ be
its  decomposition. \\
        1) We have \be~\label{ageqg} dA_t\geq g(Z_t)\,dt, P \ a.s.\ee
        2) Suppose that for some $\xi \in
L^{\infty}(\f_T)$ there is a probability measure $Q^* \sim P$ with
$U_0(\xi)=E_{Q^*}[\xi+\int_0^{T}f(q_u^*)\,du]$, then $dA_t=
g(Z_t)\,dt$ and \be U_t(\xi)=U_0(\xi)+\int_0^{ t}g(Z_u)\,du-\int_0^{
t}Z_udB_u.\ee \et

 \emph {Proof.} 1) For $\xi \in L^{\infty}(\f_T)$ and any $Q\sim P$, it
follows from the decomposition that
\begin{eqnarray}
dU_t(\xi)+f(q_t)\,dt&=&
dA_t-Z_tdB_t+f(q_t)\,dt \\
&=& ~\label{qsubmart}dA_t-Z_tq_tdt+f(q_t)\,dt-Z_tdB_{t}^Q,
\end{eqnarray}
where $B^Q$ is a $Q-$Brownian motion. This implies that
$dA_t-Z_tq_tdt+f(q_t)\,dt$ defines a non-negative measure since
$U_t(\xi)+\int_0^{ \tau\wedge t} f(q_u)\,du$ is a
 $Q$-submartingale for any $Q\sim P$. Hence $$dA_t\geq
Z_tq_tdt-f(q_t)\,dt.$$
 By taking $q^n=g'(Z)1_{\{|Z|\leq n\}}$ in the above inequality and by  letting $n$ tend to infinity, we get $dA_t\geq g(Z_t)\,dt$. \\
2) If for $\xi$, there is a measure $Q^* \sim P$ with
 $U_0(\xi)=E_{Q^*}[\xi+\int_0^{T}f(q_u^*)\,du]$, then it follows from Proposition
 \ref{qmartin} that $U_t(\xi)+\int_0^{ t } f(q_u^*)\,du$ is a
 $Q^*$-martingale. Thus applying (\ref{qsubmart}) with $Q^*$, we get
 $$ dA_t=(Z_tq_t^*-f(q_t^*))\,dt \ \ \ Q^*\ a.s.$$
Since $Q^*\sim P$, we have  \be~\label{Q*}
dA_t=(Z_tq_t^*-f(q_t^*))\,dt \ \ \ P\ a.s.\ee Finally combining
(\ref{ageqg}) and (\ref{Q*}) with the definition of
 $g$, it follows that
 $$g(Z_t)\,dt\geq (Z_tq_t^*-f(q_t^*))\,dt =dA_t\geq g(Z_t)\,dt\ \ P\ a.s.$$
 \endpf

In general we can decompose $A$ further and get:

\bc ~\label{decom}For any $\xi \in L^{\infty}(\f_T)$, there exists
an increasing predictable process $\{C_t\}_{0\leq t\leq T}$ with
$C_0=0$ such that \be
U_t(\xi)=U_0(\xi)+\int_0^tg(Z_u)\,du-\int_0^tZ_udB_u+C_t.\ee \ec

Our main result is the following.
 \bt~\label{equiv} Let $U$ be the dynamic utility function
defined by (\ref{ut}). Then the following are equivalent:
\begin{enumerate}
    \item $\underline{\lim}_{|x|\rightarrow \infty}\frac{f(x)}{|x|^2}>
    0$;
    \item $\overline{\lim}_{|z|\rightarrow \infty}\frac{g(z)}{|z|^2}<
    \infty$;
    \item For all $k> 0$, the set $\{Q\ |C_0(Q)\leq k \}$ is weakly
    compact;
    \item For all $\xi \in L^{\infty}(\f_T)$, there exists a measure
    $Q \ll P$ such that $U_0(\xi)=E_{Q}\Big[\xi+\int_0^{T}f(q_u)\,du\Big]$;
    \item For all $\xi \in L^{\infty}(\f_T)$, there exists a measure
    $Q \sim P$ such that $U_0(\xi)=E_{Q}\Big[\xi+\int_0^{T}f(q_u)\,du\Big]$;
    \item For all $\xi \in L^{\infty}(\f_T)$, the BSDE
    $dY_t=g(Z_t)\,dt-Z_tdB_t$ has a unique bounded solution with
    $Y_T=\xi$.
    \item $U_0$ is strictly monotone.

\end{enumerate}      \et

\emph {Proof.} 1 $\Leftrightarrow$ 2: Point 1 implies that there
exist positive constants $a, b \in R_+$ such that $f(x)\geq
a|x|^2-b.$ We then get$$g(z)=\sup_{x\in R^d}(zx-f(x))\leq\sup_{x\in
R^d}(zx-a|x|^2+b)\leq \frac{1}{4a}|z|^2+b$$ which shows that
$\overline{\lim}_{z\rightarrow \infty}\frac{g(z)}{|z|^2}<
    \infty$. The proof of the implication $2 \Rightarrow 1$ is similar.

    1 $\Rightarrow$ 3:   It suffices to verify that for any $k>0$, $\Big\{\frac{dQ}{dP}\ \Big|C_0(Q)=E_Q\Big[\int_0^Tf(q_u)\,du\Big]\leq k\Big\}$ is uniformly integrable.
The Dunford-Pettis theorem then shows that the set is weakly
compact.

 Since $f(x)\geq a|x|^2-b,$ we get
 $$k \geq E_Q\left[\int_0^Tf(q_u)\,du\right]\geq aE_Q\left[\int_0^T|q_u|^2\,du\right]-b.$$
 Therefore,$$\frac{1}{2}E_Q\left[\int_0^T|q_u|^2\,du\right]\leq \alpha,$$
 where $\alpha=\frac{k+b}{2a}$ is a positive constant independent of $Q$.
It follows from \bea*
\frac{1}{2}E_Q\left[\int_0^T|q_u|^2\,du\right]&=&E_Q\left[\int_0^Tq_udB_u^Q+\frac{1}{2}\int_0^T|q_u|^2\,du\right]\\
&=& E_Q\left[\int_0^Tq_udB_u-\frac{1}{2}\int_0^T|q_u|^2\,du\right] \\
&=& E_Q\left[\log \frac{dQ}{dP}\right]\eea* that for any $k>0$,
\be\bigg\{\frac{dQ}{dP}\;\bigg|E_Q\left[\int_0^Tf(q_u)\,du\right]\leq
k\bigg\} \subseteq
\bigg\{\frac{dQ}{dP}\;\bigg|E_P\left[\frac{dQ}{dP}\log
\frac{dQ}{dP}\right]\leq \alpha\bigg\}.\ee From the de la Vall\'{e}e
Poussin theorem, we conclude that
$$\bigg\{\frac{dQ}{dP}\ \bigg|E_Q\left[\int_0^Tf(q_u)\,du\right]\leq k\bigg\} \mbox{ is
uniformly integrable}.$$

3 $\Rightarrow$ 1 We prove it by the contradiction. Suppose
$\underline{\lim}_{|x|\rightarrow \infty}\frac{f(x)}{|x|^2}=
    0$, then there exists a sequence $\{x_n\}_{n=0}^{\infty}$ such
    that $\lim_{n\rightarrow\infty}|x_n|=\infty$ and $\lim_{n\rightarrow\infty}\frac{f(x_n)}{|x_n|^2}=0$. Put
    $q_n=x_n1_{[0,\delta_n\wedge T]}$ where $\delta_n=1\Big/\left(\sqrt{\frac{f(x_n)}{|x_n|^2}}|x_n|^2\right).$
    It follows from
    \be C_0(Q_n)=E_{Q_n}\left[\int_0^Tf(q_n(u))\,du\right]\le \sqrt{\frac{f(x_n)}{|x_n|^2}}\rightarrow 0,\ee
     that for all $k>0$, there exists $N>0$ such that the sequence $\{\frac{dQ_n}{dP}\}_{n=N}^{\infty}\subseteq \{\frac{dQ}{dP}\ |C_0(Q)\leq k \}$.
Furthermore, we have
\be\int_0^T|q_n|^2(u)\,du=\left(1\bigg/\sqrt{\frac{f(x_n)}{|x_n|^2}}\right)
\wedge \left(x_n^2 T\right) \rightarrow
 \infty,\ee which shows that $\frac{dQ_n}{dP}=\mathcal{E}(q_n\cdot B)_T\rightarrow
 0$, a.s. as $n\rightarrow
 \infty$.
Thus $\{\frac{dQ_n}{dP}\}_{n=N}^{\infty}$   is not
 uniformly integrable.

3 $\Leftrightarrow$ 4: It is a conclusion induced by the James'
theorem as shown in Jouini-Schachermayer-Touzi's work \cite{JST}.

    4 $\Leftrightarrow$ 5: It is obvious that point 5 implies point 4. For the proof of the inverse implication, we use the fact that condition 4 is equivalent to condition 2.
    In this case, by convexity, there exists a positive constant $c$ such that $|g'(z)|\leq c(|z|+1)$. For any $\xi \in
    L^{\infty}(\f_T)$, there is a measure $Q\ll P$ such that
    $U_0(\xi)=E_{Q}[\xi+\int_0^{T}f(q_u)\,du]$, then, by Proposition ~\ref{qmartin}, $U_t(\xi)+\int_0^{\tau\wedge t } f(q_u)\,du$ is a
 $Q$-martingale  where $\tau=\inf\{t\in
[0,T]\ |\ \mathcal{E}(q\cdot B)_t=0\}\wedge T$. It follows from
    (\ref{qsubmart}) that
$$dA_t=(Z_tq_t-f(q_t))\,dt \ \ \ m\otimes Q\ a.s.\mbox{ on }[0,\tau],$$
where $m$ is the Lebesgue measure on $[0,T]$.  Since $dA_t\geq
g(Z_t)\,dt$, $m\otimes Q$\ a.s., we get $$g(Z_t)=Z_tq_t-f(q_t)\ \
m\otimes Q\ a.s.,$$ which implies $q_t=g'(Z_t)$ on $[0,\tau]$. We
then have
    \bea*\int_0^{\tau}|q_u|^2\,du&=&\int_0^{\tau}(g'(Z_u))^2\,du\\
    &\leq& c^2\int_0^{\tau}(1+|Z_u|)^2\,du< \infty,
    \eea*
which means
$P\left\{\frac{dQ}{dP}=0\right\}=P\Big\{\int_0^{\tau}|q_u|^2\,du=\infty\Big\}=0$.
Hence $Q\sim P$.

5 $\Rightarrow$ 6: For a given $\xi \in
    L^{\infty}(\f_T)$, if there exists a measure
    $Q \sim P$ such that $U_0(\xi)=E_{Q}\left[\xi+\int_0^{T}f(q_u)\,du\right]$,
    it follows from Lemma \ref{atgz} that $\left\{U_t,Z_t\right\}_{0\leq t\leq T}$ is a
    solution of the following BSDE:

\be~\label{BSDE} \left\{\ba{l}dY_t=g(z_t)\,dt-z_tdB_t,\q\q 0\leq t\leq T;\\
 Y_T=\xi,\q \xi \in L^{\infty}(\f_T); \\ Y \hbox{ \rm is bounded }
    ,\ea\right.
\ee where $E\left[\int_0^T|z_t|^2\,dt\right]<\infty$ and
$E\left[\int_0^Tg(z_t)\,dt\right]<\infty$. Since, as we have proved
above, condition 5 implies $\overline{\lim}_{z\rightarrow
\infty}\frac{g(z)}{|z|^2}<
    \infty$, the BSDE has a unique bounded solution
    according to Kobylanski \cite{Kobylanski}.

6 $\Rightarrow$ 2 We will prove this in the next section. See
Theorem \ref{nonexistence}.

5 $\Rightarrow$ 7 For any $\xi \in L^{\infty}(\f_T)$, there exists
an equivalent measure $Q\sim P$ such that
$U_0(\xi)=E_{Q}[\xi+\int_0^{T}f(q_u)\,du]$ with
$\frac{dQ}{dP}=\mathcal{E}(q\cdot B)$.

 Suppose that
$U_0(\eta)=U_0(\xi)$ for some $\eta \in L^{\infty}(\f_T)$ with
$\eta\leq \xi$, $P$ a.s. Since
$$U_0(\eta)\leq E_{Q}\left[\eta+\int_0^{T}f(q_u)\,du\right]\leq E_{Q}\left[\xi+\int_0^{T}f(q_u)\,du\right]=U_0(\xi),$$
we have $E_Q[\xi-\eta]=0$, hence $\xi=\eta$, $Q$ a.s. Thus
$\xi=\eta$, $P$ a.s. and $U_0$ is strictly monotone.

7 $\Rightarrow$ 2 See Remark ~\ref{not-strict}, Remark ~\ref{minex}
or Example ~\ref{mincont}.
\endpf

We have proved that in the case when the generator $g$ is at most
quadratic, the dynamic utility function $U$ is the solution of BSDE
(\ref{BSDE}). In general, however, we have the following inequality.

\bl ~\label{ugeqy}For any $\xi \in L^{\infty}(\f_T)$, if BSDE
(\ref{BSDE}) has a bounded solution $Y$, then we have $U(\xi)\geq
Y$.\el

\emph{Proof.} $Y$ is bounded. The following calculation is therefore
justified: \bea* E_Q\left[\xi+\int_t^Tf(q_u)\,du\Big|\f_t\right]&=&
Y_t+E_Q\left[\int_t^Tg(Z_u)\,du-\int_t^TZ_udB_u+\int_t^Tf(q_u)\,du\Big|\f_t\right]\\
&=&Y_t+E_Q\left[\int_t^T[g(Z_u)-Z_uq_u+f(q_u)]\,du\Big|\f_t\right]\\
&\geq& Y_t, \mbox{ for any } Q\sim P \mbox{ with }
E_Q\left[\int_0^Tf(q_u)\,du\right]<\infty.\eea* \endpf
\section{Backward SDEs with superquadratic growth.}
In this section, we discuss the following BSDE($g$, $\xi$):
\be~\label{BSDEgsq} \left\{\ba{l}dY_t=g(Z_t)\,dt-Z_tdB_t;\\
 Y_T=\xi,\q \xi \in L^{\infty}(\f_T),\ea\right.
\ee where  $g:R^{ d}\rightarrow R_{+}\cup\{+\infty\}$ is convex with
$g(0)=0$ and superquadratic $\overline{\lim}_{|z|\rightarrow
\infty}\frac{g(z)}{|z|^2}=
    \infty.$ A pair of predictable processes ($Y$, $Z$) is called a bounded solution to BSDE (\ref{BSDEgsq}) if \\$Y:\Omega\times[0,T]\rightarrow R$ is bounded and
\\$Z:\Omega\times[0,T]\rightarrow R^{ d} $ is such that $E\left[\int_0^Tg(Z_t)\,dt\right]<\infty$.

Here for simplicity, we consider the BSDE with $d=1$. However, the
results remain valid for $d>1$.

\subsection{Non-existence of the solution}
Different from the BSDEs with at most quadratic growth, the solution
to the BSDE with super-quadratic growth does not always exist.
 \bt{\bf (Non-existence)}~\label{nonexistence} There
exists $\eta\in L^{\infty}(\f_T)$ such that BSDE(\ref{BSDEgsq}) with
superquadratic growth has no bounded solution. \et

\emph{Proof.} The proof is divided into 4 steps.\\
{\bf Step 1.} We construct a pair of processes $(X,Z)$, a measure
$Q$ as well as a bounded  random variable $\xi$.

 Since $\overline{\lim}_{|z|\rightarrow
\infty}\frac{g(z)}{|z|^2}=
    \infty$, there exists a sequence $\{z_k\}_{k=1}^{\infty}$ such
    that $\lim_{k\rightarrow \infty}|z_k|=\infty$ and $g(z_k)\geq k|z_k|^2.$ Without loss of generality, we suppose  $z_k>0$. The other case is left to the reader. Thus we have \be ~\label{g'} g'(z_k)\geq \frac{g(z_k)}{z_k}\geq k\cdot
z_k.\ee
    We put $Z_u\triangleq \sum_{n=1}^{\infty}z_n 1_{[\sum_{k<n}\delta_k,\sum_{k\leq n}\delta_k
    )}(u)$ where $$\delta_k=\frac{1}{\alpha z_kg'(z_k)k^2} $$ and we
set
    $$\alpha=\sum_{k=1}^{\infty}\frac{1}{Tz_kg'(z_k)k^2}<\infty,$$
in order to have
$$\sum_{k\geq 1}\d_k=\sum_{k\geq 1}\frac{1}{\alpha z_kg'(z_k)k^2}=T.$$
    Then from (\ref{g'}), we have
        $$\int_0^Tg(Z_u)\,du=\sum_{k\geq 1}g(z_k)\d_k\leq \sum_{k \geq 1} z_kg'(z_k)\d_k=\sum_{k\geq 1}\frac{1}{\alpha k^2}<\infty,$$
    $$\int_0^T|Z_u|^2\,du=\sum_{k\geq 1}|z_k|^2\d_k\leq\sum_{k\geq 1}z_kg'(z_k)\d_k\frac{1}{k}= \sum_{k\geq 1}\frac{1}{\alpha k^3}<\infty.$$

Let $q_t=g'(Z_t)$. It follows from
$$\int_0^T|q_u|^2\,du=\sum_{k\geq 1}(g'(z_k))^2\d_k\geq \sum_{k\geq 1}kg'(z_k)z_k\d_k=\sum_{k\geq 1}\frac{1}{\alpha k}=+\infty,$$
that $\lim_{t\rightarrow T}\mathcal{E}(q\cdot B)_t=0$ and
$\mathcal{E}(q\cdot B)_t>0,\ P\hbox{\ a.s. for any } t<T.$

Let $X_t=\int_0^tg(Z_u)\,du-\int_0^tZ_udB_u.$ We stop $X$ at a
random time $\si$  \be~\label{defsi}\si\triangleq \inf\{t\in
[0,T]\mid \mathcal{E}(q\cdot B)_t\geq n\}\wedge \inf\{t\in
[0,T]\mid|X_t|\geq n\}\wedge T\ee where $n$ is a positive constant
which is sufficiently large to ensure that $P(\si=T)>0$. We then set
a measure $Q^*$ with
$E_P\left[\frac{dQ^*}{dP}\Big|\f_t\right]=\mathcal{E}(q^*\cdot
B)_{t}$ and $q^*_t=g'(Z_t)1_{\{t\leq \si\}}.$\\ We define
$\xi=X_{\si}\in L^{\infty}(\f_T).$

 {\bf Step 2.} The measure $Q^*\ll P$ but it is not
 equivalent to $P$.

  Set $A_1=\{\si=T\}$. Then
 $$Q^*(A_1)=\int_{A_1}\mathcal{E}(q^{}\cdot B)_{\si}dP=\int_{A_1}\mathcal{E}(q^{}\cdot B)_{T}dP=0$$
 while $P(A_1)>0$.
Thus we have $Q^*\nsim P$ and $Q^*\ll P$. However, $Q^*\bigotimes
m\sim P\bigotimes m$ where $m$ is the Lebesgue measure since
$Q^*\sim P$ on $\f_t$ for all $  t<T$. Clearly $(X^{\si}_t,\
Z_t1_{\{t\leq \si\}})_{0\le t\le T}$ is a bounded solution of BSDE
($g$, $\xi$) where $X^{\si}_t=X_{\si \wedge t}$.

{\bf Step 3.} In this step we prove that the dynamic utility
function $U(\xi)$ is the bounded solution of BSDE ($g$,
$\xi$) (\ref{BSDEgsq}) and
$U_t(\xi)=E_{Q^*}\left[\xi+\int_t^Tf(q_u^*)\,du\Big|\f_t\right]$ for
any $t<T$.

As $X^{\si}$ is a bounded solution of BSDE ($g$, $\xi$), we get
  \begin{eqnarray}~\label{uleqy} U_t(\xi)&\leq&
E_{Q^*}\Big[\xi+\int_t^Tf(q_u^*)\,du\Big|\f_t\Big]\nonumber \\
&=&E_{Q^*}\Big[\xi+\int_{t\wedge \si}^{\si}f(q_u^{*})\,du\Big|\f_t\Big]\nonumber \\
&=& X_t^{\si}+E_{Q^*}\Big[\int_{t\wedge \si}^{\si}(f(q_u^*)+g(Z_u))\,du-\int_{t\wedge \si}^{\si}Z_udB_u\Big|\f_t\Big]\nonumber\\
&=& X_t^{\si}+E_{Q^*}\Big[\int_{t\wedge \si}^{\si}[f(q_u^*)+g(Z_u)-Z_uq_u^*]\,du\Big|\f_t\Big]\nonumber\\
&=& X_t^{\si},\ \ Q^*\ a.s.\end{eqnarray} hence $P$ a.s. because
$Q^*\sim P$ on $\f_t$ for $t<T$. Combining Lemma \ref{ugeqy} with
inequality (\ref{uleqy}), we deduce that\be~\label{ueqy} U_t(\xi)=
E_{Q^*}\Big[\xi+\int_t^Tf(q_u^*)\,du\Big|\f_t\Big]= X_t^{\si},\ \ P
\hbox{ a.s. for all }t\in[0,T).\ee Set $\eta=\xi+h$ where $h \in
L_{+}^{\infty}(\f_T)$, $P[h>0]>0$ and $h\cdot\mathcal{E}(q^*\cdot
B)_{\si}=0$.

{\bf Step 4.} We show that $U_t(\xi)=U_t(\eta),\ P$ a.s. for any
$t<T$ and hence BSDE ($g$,$\eta$) has no solution.

 It follows from $\eta=\xi,\ Q^*$-a.s. that\bea*
U_t(\eta)&\leq&
E_{Q^*}\Big[\eta+\int_t^Tf(q_u^*)\,du\Big|\f_t\Big]\\
&=& E_{Q^*}\Big[\xi+\int_t^Tf(q_u^*)\,du\Big|\f_t\Big]\\
&=&U_t(\xi)\ \ \hbox{ for any }t< T.\eea*

Notice that $U$ is monotone, i.e., $U_t(\xi)\leq U_t(\eta)$, and so
we have $U_t(\xi)=U_t(\eta),\ P$ a.s. for any $t<T$.

Suppose $Y$ is a bounded solution of BSDE ($g$, $\eta$), then  we
have for $t<T$,
$$X_t^{\si}=U_t(\xi)=U_t(\eta)\geq Y_t,$$
and hence
$$\eta=Y_T=\lim_{t\rightarrow T}Y_t\leq \lim_{t\rightarrow
T}X_t^{\si}=X_T^{\si}=\xi, \ P\ a.s.,$$ a contradiction to the fact
that $P[\eta>\xi]>0$. Therefore, BSDE ($g$, $\eta$) has no solution.
\endpf

\br  From this theorem, together with what we have proved in Theorem
\ref{equiv} we get that BSDE ($g$, $\xi$) has a solution for all
$\xi\in L^{\infty}(\f_T)$ if and only if $g$ is at most
quadratic.\er

\br ~\label{not-strict} From the proof, we get $\eta\geq\xi$ with
$P(\eta>\xi)>0$ and $U_0(\xi)=U_0(\eta)$. Thus the utility function
$U_0$ is NOT strictly monotone when
$\underline{\lim}_{|x|\rightarrow \infty}\frac{f(x)}{|x|^2}=0$.\er

Although the BSDE ($g$, $\xi$) (\ref{BSDEgsq}) does not always have a
solution, in the following case it has.

\bde We say that a random variable $\xi\in L^{\infty}(\f_T)$ is
minimal if $\eta \leq \xi$ and $P[\eta<\xi]>0$ imply
$U_0(\eta)<U_0(\xi)$. \ede

\bt Let $\xi\in L^{\infty}(\f_T)$ be minimal. Then $U(\xi)$ is a
solution of BSDE ($g$, $\xi$). \et

\emph{Proof.} We prove it by contradiction. Let $\xi\in
L^{\infty}(\f_T)$ be minimal and suppose $U(\xi)$ is not a solution
of BSDE ($g$, $\xi$). Then it follows from Corollary \ref{decom}
that there exists an increasing process $C$ with $C_0=0$ such that
$P[C_T>0]>0$ and \be
U_t(\xi)=\xi-\int_t^Tg(Z_u)\,du+\int_t^TZ_udB_u-C_T+C_t. \ee Define
$\tau:=\inf\{t\in [0,T]| C_t\ge k\}\wedge T$, where $k>0$ is such
that $P[C_\tau>0]>0$. Since $C$ may have jumps, $C_{\tau}$ can be
unbounded. However, $\tau$ is predictable so there exists
$\{\tau_n\}_{n=1}^{\infty}$ such that $\tau_n\uparrow \tau$ and
$\tau_n<\tau$ on $\{\tau>0\}$. It follows that $C_{\tau_n}\leq k$
and $P[C_{\tau_n}>0]>0$ for $n$ big enough. Denote by $\sigma$ a
stopping time $\tau_n$ for $n$ big enough, then we have
$$
U_t(\xi)-C_t=U_\sigma(\xi)-C_\sigma-\int_t^\sigma g(Z_u)\,du+\int_t^\sigma Z_udB_u,
$$
which implies that $(U_{t\wedge\sigma}(\xi)-C_{t\wedge\sigma}, Z_t1_{\{t\le\sigma\}})_{0\le t\le T}$ is a solution of
BSDE ($g$, $U_\sigma(\xi)-C_\sigma$). Thus by Lemma \ref{ugeqy}, we deduce
$$U_0(\xi)=U_0(\xi)-C_0\leq U_0(U_\sigma(\xi)-C_\sigma).$$

On the other hand, it is clear that $U_0(\xi)\geq U_0(U_\sigma(\xi)-C_\sigma)$.
 Therefore, we have $$U_0(\xi)=U_0(U_\sigma(\xi)-C_\sigma).$$
It follows from the above equality, the translability and the
time-consistency  of the dynamic utility function that
$$U_0(\xi)=U_0\left(U_\sigma(\xi)-C_\sigma\right)=U_0\left(U_\sigma(\xi-C_\sigma)\right)=U_0(\xi-C_\sigma).$$
This is a contradiction to the fact that $\xi$ is minimal. \endpf

\br For $g$ with at most quadratic growth
$\overline{\lim}_{|z|\rightarrow \infty}\frac{g(z)}{|z|^2}<
    \infty$, it follows from Theorem \ref{equiv} that $\xi$ is minimal for all $\xi \in L^{\infty}(\f_T)$.
\er
 If $g$ is superquadratic, there exists a bounded random variable $\zeta$
such that $U(\zeta)$ is a solution of BSDE ($g$, $\zeta$) and
$\zeta$ is not minimal. See Example ~\ref{mincont}.

\subsection{Non-uniqueness of the Solution}
 In this subsection, we shall prove that if the BSDE has a bounded
solution, the bounded solution is not unique. The main reason is
that the generator $g$ is superquadratic which makes
$\int_0^{t}g(Z_r)dr$ grow much faster than $\int_0^{t}Z_rdB_r$.
Following this observation, we can construct other solutions.

\bt~\label{non-uniqueness}{\bf (Non-uniqueness)} If the BSDE ($g$,
$\xi$) with superquadratic growth has a bounded solution $Y$ for a
$\xi\in L^{\infty}(\f_T)$, then for each $y<Y_0$, there are
infinitely many bounded solutions $\{X_t\}_{0\leq t\leq T}$ with
$X_0=y$.\et

\emph{Proof.}  Suppose $(Y,Z)$ is a bounded solution of BSDE ($g$,
$\xi$). Divide the time interval $[0,T]$ into
$[T(1-2^{-n}),T(1-2^{-n-1}))$, where $n=0,1,2,...$  and denote
$\alpha_n=T(1-2^{-n})$. Suppose the new solution $(X,Z')$ has been
constructed on $[0,\alpha_n]$ with $X_0=y$ where $y<Y_0$ such that
$X_{\a_n}\leq Y_{\a_n}$ $P$ a.s. Let us construct $(X,Z')$ on the
time interval $[\alpha_n,\alpha_{n+1})$.

Our idea is the following. Since $g$ is superquadratic, we can
construct a process $X_{\a_n}+V_t,\ t\in [\alpha_n,\alpha_{n+1})$
such that $\lim_{t\rightarrow \a_{n+1}}V_t=+\infty$, $P$ a.s. and
for any $0<\e<1$, $V_t$ exceeds downwards $-2^{-n-1}\e$ with a very
small probability. The fact that the solution $Y$ is bounded implies
that it is touched by the process $X_{\a_n}+V_t$ because
$X_{\a_n}\leq Y_{\a_n}$. We then get a new solution $X_t$ on this
time interval $[\a_n, \a_{n+1}]$ by stopping $X_{\a_n}+V_t$ when it
reaches $Y$.

First, let us construct the process $V_t$.

 It follows from $\overline{\lim}_{z\rightarrow
\infty}\frac{g(z)}{|z|^2}=
    \infty$ that there exists a sequence $\{x_k\}_{k=0}^{\infty}$
    such that for any $k\geq 0$,
    \begin{enumerate}
        \item  $g(x_k)\geq 4^n x_k^2$;
        \item
        $x_k^2\geq\frac{1}{(\theta^k-\theta^{k+1})\theta^k\delta_n}$
        where $\theta\in (0,1)$ is a constant and $\delta_n=\alpha_{n+1}-\alpha_n=2^{-n-1}T.$
    \end{enumerate}

Set $b_t=\sum_{k=0}^{\infty} x_k1_{[\alpha_{n+1}-\th^k\d_n, \
\alpha_{n+1}-\th^{k+1}\d_n)}(t)$ and
$V_t=\int_{\a_n}^tg(b_u)\,du-\int_{\a_n}^tb_udB_u$ for any $t\in
[\a_n,\  \a_{n+1})$.

We then have for $t\in[\alpha_{n+1}-\th^{N+1}\d_n, \
\alpha_{n+1}-\th^{N+2}\d_n)$,

\be\int_{\a_n}^tb_u^2\,du\geq
\sum_{k=0}^Nx_k^2(\th^k-\th^{k+1})\d_n\geq\sum_{k=0}^N\frac{1}{\th^k},\ee
\be\int_{\a_n}^tg(b_u)\,du\geq4^n
\int_{\a_n}^tb_u^2\,du\geq4^n\sum_{k=0}^N\frac{1}{\th^k}.\ee Thus
$\lim_{t\rightarrow \a_{n+1}}\int_{\a_n}^tb_u^2\,du=\infty$ and
$\lim_{t\rightarrow \a_{n+1}}\int_{\a_n}^tg(b_u)\,du= \infty$.

{\bf Step 1.} \emph{We have $\lim_{t\rightarrow
\a_{n+1}}V_t=+\infty$ $P$ a.s.}

 Define $\phi(t)=\int_{\a_n}^tb_u^2\,du$
for $t\in [\a_n, \a_{n+1})$. Then $\phi$ is strictly increasing with
$\phi(\a_n)=0$ and
$$\lim_{t\rightarrow \a_{n+1}}\phi(t)=+\infty.$$ Setting
$B_t^{*}\triangleq \int_{\a_n}^{\phi^{-1}(t)}b_udB_u$, we get a time
changed Brownian motion with respect to the filtration
$\{\f^{B^{*}}\}$. It follows from the construction of $V$ that
$$V_t\geq 4^n\phi(t)-B^{*}_{\phi(t)}=\phi(t)\bigg[4^n-\frac{B^{*}_{\phi(t)}}{\phi(t)}\bigg],$$
which implies that \be\lim_{t\rightarrow \a_{n+1}} V_t= +\infty,\ P\
\mbox{ a.s.}\ee since $$\lim_{t\rightarrow \a_{n+1}}
\frac{B^{*}_{\phi(t)}}{\phi(t)}= 0,\ P\ \mbox{ a.s.}$$


 Now we
estimate the probability that $V_t$ reaches a small negative number
$-2^{-n-1}\e$.

{\bf Step 2.} \emph{Calculate the probability$$P(\{\o\in \Omega \mid
\exists \,t\in [\a_n, \a_{n+1}) \hbox{ such that }
V_t(\o)<-2^{-n-1}\e\}).$$} Applying the submartingale inequality, we
deduce that
\begin{eqnarray}~\label{exceedProb}
&&P(\{\o\in \Omega\mid \exists\, t\in [\a_n, \a_{n+1}) \hbox{ such
that } V_t(\o)<-2^{-n-1}\e\}) \nonumber \\
&=& P(\{\o\in \Omega\mid\exists\, t\in [\a_n, \a_{n+1}) \hbox{ such
that } 4^n\phi(t)-B_{\phi(t)}^{*}<-2^{-n-1}\e\}) \nonumber \\
&=& P(\{\o\in \Omega\mid \exists \,s\in[0,\infty) \hbox{ such
that } 4^n s-B_{s}^{*}<-2^{-n-1}\e\}) \nonumber \\
&\leq &\exp\{-2^n\e\}.
\end{eqnarray}

{\bf Step 3.} \emph{Construct the new solution $(X_t,Z'_t)$ for all
$t\in [\a_n, \a_{n+1}]$.}

Define $$\tau_1\triangleq \inf\{t\geq
\a_n\,|\,V_t=-2^{-n-1}\e\}\wedge \a_{n+1}$$ and $$\tau_2\triangleq
\inf\{t\geq \a_n\,|\, X_{\a_n}+V_t\geq Y_t\}\wedge \a_{n+1}$$ which
are the stopping times when the process $X_{\a_n}+V_t$ touches
$X_{\a_n}-2^{-n-1}\e$ and $Y_t$ respectively. It follows from
$\lim_{t\rightarrow \a_{n+1}} V_t= +\infty\ P\ \mbox{ a.s.}$ that
$P[\tau_2<\a_{n+1}]=1$. Define
$$\tau_3\triangleq \inf\{t\geq \tau_1\,|\,X_{\a_n}-2^{-n-1}\e=
Y_t\}\wedge \a_{n+1}.$$

Now we have three cases \be\left\{
  \begin{array}{ll}
    \tau_1<\tau_2, \tau_3<\a_{n+1}, & \hbox{put }Z'_t(\o)=b_t1_{\{t\leq \tau_1\}}+Z_t1_{\{t> \tau_3\}} \hbox{;}  \\
    \tau_1<\tau_2, \tau_3=\a_{n+1}, &  \hbox{put }Z'_t(\o)=b_t1_{\{t\leq \tau_1\}}\hbox{;} \\
    \tau_1\geq\tau_2,               & \hbox{put }Z'_t(\o)=b_t1_{\{t\leq \tau_2\}}+Z_t1_{\{t> \tau_2\}}\hbox{,}
  \end{array}
\right.\ee where $(Y,Z)$ is the original bounded solution of the
BSDE ($g$, $\xi$).

 Thus we
get\begin{eqnarray}Z'_t&=&1_{\{\tau_1<\tau_2,
\tau_3<\a_{n+1}\}}(b_t1_{\{t\leq \tau_1\}}+Z_t1_{\{t>
\tau_3\}})\nonumber \\
 &+&1_{\{\tau_1<\tau_2, \tau_3=\a_{n+1}\}}b_t1_{\{t\leq
 \tau_1\}}\nonumber \\
 &+&1_{\{\tau_1\geq\tau_2\}}(b_t1_{\{t\leq
\tau_2\}}+Z_t1_{\{t> \tau_2\}})\\
&=& 1_{\{t\leq \tau_1\wedge
\tau_2\}}b_t+[1_{\{\tau_1<\tau_2,\tau_3<\a_{n+1},t>\tau_3\}}+1_{\{\tau_1\geq\tau_2,t>\tau_2\}}]Z_t.
\end{eqnarray}
Obviously, $Z'$ is a predictable process.

Set \be X_t\triangleq
X_{\a_n}+\int_{\a_n}^tg(Z'_u)\,du-\int_{\a_n}^tZ'_udB_u \ee for all
$t\in [\a_n, \a_{n+1}]$.

{\bf Step 4.} \emph{Some properties of $X$.}

It follows from the construction that $\{X_t\}_{\a_n\leq t\leq
\a_{n+1}}$ has the following properties:
\begin{enumerate}
    \item  $X_t1_{\{\tau_2\leq \tau_1,t\geq \tau_2\}}=Y_t1_{\{\tau_2\leq \tau_1,t\geq
    \tau_2\}}$;
    \item  $X_t1_{\{\tau_1<\tau_2,\tau_3\geq t\geq \tau_1\}}=(X_{\a_{n}}-2^{-n-1}\e)1_{\{\tau_1<\tau_2,\tau_3\geq t\geq
    \tau_1\}}$;
    \item  $
    X_t1_{\{\tau_1<\tau_2,t>\tau_3\}}=Y_t1_{\{\tau_1<\tau_2,t>\tau_3\}}$.
\end{enumerate}

Therefore, we have
\begin{eqnarray}X_{\a_{n+1}}&=&Y_{\a_{n+1}}(1_{\{\tau_2\leq\tau_1\}}+1_{\{\tau_1<\tau_2,\tau_3<\a_{n+1}\}})\nonumber \\
&+&(X_{\a_n}-2^{-n-1}\e)1_{\{\tau_1<\tau_2,\tau_3=\a_{n+1}\}}.
\end{eqnarray}
So the induction assumption $X_{\a_n}\leq Y_{\a_n}$ implies
$X_{\a_{n+1}}\leq Y_{\a_{n+1}}$. It is also clear that the new
solution $X$ is bounded by $\parallel Y\parallel_{\infty}+|y|+\e$.

 Set $A_n\triangleq\{\o\in\Omega\mid\tau_1<\tau_2,
\tau_3=\a_{n+1}\}$. Then $P(A_n)$ is the probability that
$X_{\a_{n+1}}$ is not equal to $Y_{\a_{n+1}}$. From
(\ref{exceedProb}), we get
$$P(A_n)\leq P(\{\o\in \Omega\mid \exists\, t\in [\a_n, \a_{n+1}) \hbox{ such
that } V_t(\o)<-2^{-n-1}\e\})\leq \exp\{-2^n\e\}.$$

Since $\sum_{n=0}^{\infty}\exp\{-2^n\e\}<+\infty$, the
Borel-Cantelli Lemma implies that \be
P(\cap_{n=0}^{\infty}\cup_{k\geq n}A_k)=0,\ee which shows
$X_T=Y_T=\xi,\q P$ a.s.


To sum up, $(X,Z')$ is indeed a new bounded solution with $X_0=y$.

The construction used many different constants. It is clear that
this yields infinitely many different solutions.
\endpf

Notice that in the proof we only use the fact that $g$ is
superquadratic to guarantee that the new solution $X$ is bounded
below. This shows if g is at least quadratic , i.e.
$\overline{\lim}_{|z|\rightarrow \infty}\frac{g(z)}{|z|^2}>0$, we
can construct a process $V_t$ such that $\lim_{t\rightarrow
\a_{n+1}}V_t=+\infty$ as well. Thus we have the following
conclusion. \bc Suppose $g$ is at least quadratic
$\overline{\lim}_{|z|\rightarrow \infty}\frac{g(z)}{|z|^2}>0$ and,
for $\xi \in L^{\infty}(\f_T)$, $Y$ is a bounded solution of the
BSDE ($g$, $\xi$), then for each $y<Y_0$, there exists infinitely
many solutions $X$ which are {\rm \bf bounded above} with
$X_0=y$.\ec

\subsection{Non-stability of the solutions}

The monotone stability plays an important role in the study of
quadratic BSDEs (See, e.g., \cite{Kobylanski,BH1}). Here we shall
show that the same type of monotone stability does not hold.

\bt {\bf (Non-stability)} Suppose $\overline{\lim}_{z\rightarrow
\infty}\frac{g(z)}{|z|^2}=
    \infty$. Then there exists a sequence of solutions $\{Y^k\}_{k=1}^{\infty}$ of  BSDEs ($g$, $\xi_k$) which increasingly
    and boundedly converges to $Y$  such that $Y$ is not a solution of BSDE ($g$, $\xi$), where $\xi$ is the $L^{\infty}$ limit of $\{\xi_k\}_{k=1}^{\infty}$.\et

\emph{Proof. } It follows from $\overline{\lim}_{z\rightarrow
\infty}\frac{g(z)}{|z|^2}=
    \infty$ that there exists a sequence
$\{z_k\}_{k=1}^{\infty}$ with $|z_k|\rightarrow +\infty$ such that $g(z_k)\geq
\max\{16^kT|z_k|^2,2^{k+1}T\}.$ W.l.o.g., we suppose that $z_k>0$.

     Denote $\alpha_k:=\lceil g(z_k)\rceil$ where $\lceil\cdot\rceil$ is the ceiling function and put $Z^k(t)\triangleq z_k1_{\cup_{i=1}^{\alpha_k}[\frac{T}{\alpha_k}i-\frac{T}{\alpha_k^2},
     \frac{T}{\alpha_k}i]}(t)$ for all $0\leq t\leq T$. Then it follows
     that
     \bea* E\left[\int_0^t{Z}_u^kdB_u\right]^2&\leq& \int_0^T({Z}_u^k)^2\,du\\ &=&
(z_k)^2\frac{T}{\alpha_k}\leq 16^{-k}\rightarrow 0,\  \mbox{for }
k\rightarrow \infty.
     \eea*
 However, we have \be~\label{estimgz} \int_0^tg(Z_u^k)\,du =
\int_0^tg(z_k)1_{\cup_{i=1}^{\alpha_k}[\frac{T}{\alpha_k}i-\frac{T}{\alpha_k^2},
     \frac{T}{\alpha_k}i]}(u)\,du \in \Big[\frac{g(z_k)}{\alpha_k}\left(t-\left(\frac{T}{\alpha_k}-\frac{T}{\alpha_k^2}\right)\right),\frac{g(z_k)}{\alpha_k}t\Big],\ee
which implies \be\sup_{0\leq t\leq T}\Big|\int_0^tg(Z_u^k)\,du
-t\Big|\le 2^{-k}\rightarrow 0,\ee as $k
\rightarrow \infty$.


Define stopping times
\be\nu_k\triangleq\inf\Big\{t\ge 0\,\Big|\,\big|\int_0^tZ_u^kdB_u\big|>
2^{-k}\Big\}\wedge T\ee and \be~\label{defnu} \nu=\inf_{k\geq
1}\nu_k.\ee

Applying the submartingale inequality, we get \bea* P[\nu_k<T]&=&
P\left[\sup_{0\le t\le T}
\big|\int_0^tZ_u^kdB_u\big|> 2^{-k} \right]\\
&\leq& 4^k E\left[\Big(\int_0^TZ_u^kdB_u\Big)^2\right]\leq
4^{-k}.\eea*

 Thus we get
\bea*P[\nu=T]&=& 1-P[\cup_{k\geq 1}\{\nu_k<T\}]\\
&\geq& 1-\sum_{k\geq 1}P[\nu_k<T]\\
&\ge& \frac{2}{3}>0,\eea* which is due to the selection of
sufficient large $z_k,k\geq 1$.

Since $\sum_{k\geq 1}P[\nu_k<T]<\infty$, it follows from the
Borel-Cantelli Lemma that
$$P\big[\cap_{n\geq 1}\!\!\cup_{k\geq n}\{\nu_k<T\}\big]=0,$$ which means
$P[\cup_{n\geq 1}\!\!\cap_{k\geq n}\{\nu_k=T\}]=1$. It implies that,
for almost all $\o\in \Omega$, there exists $N(\o)$ such that for
any $k>N(\o)$, $\nu_k(\o)=T.$ Thus we have $P[\nu>0]=1$.

Define $y_t^k=\int_0^{t \wedge \nu}g(Z_u^k)\,du-\int_0^{t \wedge
\nu}Z_u^kdB_u$. We then deduce that \begin{eqnarray}\label{extra} &&\sup_{0\leq t \leq T}|y_t^k-t \wedge \nu|\nonumber \\
&\leq&\sup_{0\leq t \leq T}\|\int_0^{t \wedge
\nu}g(Z_u^k)\,du-t
\wedge \nu\|+\sup_{0\leq t \leq T}\|\int_0^{t \wedge \nu}Z_u^kdB_u\|\nonumber \\
&\leq& 2\cdot 2^{-k}, \end{eqnarray}
which implies that
 \be~\label{supyk}\lim_{k\rightarrow
\infty}\sup_{0\leq t \leq T}|y_t^k-t \wedge \nu|=0.\ee

Set \be
Y_t^n=y_t^n-8+\sum_{k=1}^{n}4\cdot 2^{-(k-1)}.\ee

Notice that the stopping time $\nu_k$ is defined such that
$|\int_0^{t\wedge\nu}Z^k_udB_u|\leq 2^{-k},\forall
t\in[0,T]$. Combining (\ref{extra}) with the definition of
$\nu_k$, we get
that \bea* Y_t^{k}-Y_t^{k-1}&=& y_t^{k}-y_t^{k-1}+4\cdot 2^{-k}\\
&\geq&\[t\wedge\nu-2\cdot 2^{-k}\]-\[t\wedge\nu+2\cdot 2^{-(k-1)}\]+4\cdot 2^{-(k-1)}\\
&\ge&0\eea* which shows that $\{Y^k\}_{k=1}^{\infty}$ is a
nondecreasing sequence. Set $\xi^k=Y_T^k$ for $k\geq 1$. Then $Y^k$
is a solution of the BSDE ($g$, $\xi^k$). It follows from
(\ref{supyk}) that $Y_t^k$ converges to $t\wedge \nu$ as
$k\rightarrow\infty$ and $\xi^k\rightarrow\nu$ in $L^{\infty}$.
However, $t\wedge \nu$ is not a solution of the BSDE ($g$, $\nu$) for
$t\wedge \nu$ is an increasing process.\endpf

\br Although $t\wedge \nu$ is not the solution of the BSDE ($g$,
$\nu$), it is the dynamic utility function of $\nu $, i.e.
$U_t(\nu)=t \wedge \nu$.\er

\emph{Proof.} Indeed, setting the measure $Q^k$ such that
$E_P[\frac{dQ^k}{dP}|\f_t]=\mathcal{E}(q^k\cdot B)_t$ where
$q_t^k=g'(z_k)1_{\cup_{i=1}^{\alpha_k}[\frac{T}{\alpha_k}i-\frac{T}{\alpha_k^2},
     \frac{T}{\alpha_k}i]}(t\wedge \nu)$, we have
$$E\left[\exp\Big\{\frac{1}{2}\int_0^T(q_t^k)^2\,dt\Big\}\right]\leq \exp\Big\{\frac{1}{2}(g'(z_k))^2T\Big\}<\infty.$$
So $\mathcal{E}(q^k\cdot B)_t$ is a $P$-martingale and $Q^k$ is well
defined. Then
     \bea*
     U_t(\xi^k)&\leq&E_{Q^k}\[\xi^k+\int_t^Tf(q_u^k)\,du\|\f_t\]\\
     &=&Y_t^k+E_{Q^k}\[\int_t^{\nu}(f(q_u^k)+g(Z_u^k)-Z_u^kq_u^k)\,du\|\f_t\]\\
     &=&Y_t^k.\eea*
Thus it follows from Lemma \ref{ugeqy} that $U_t(\xi^k)=Y_t^k$. If
$k$ tends to infinity, we get $$U_t(\nu)=t\wedge \nu,$$ since
$\xi^k\rightarrow\nu$ in $L^{\infty}$.
\endpf
\br ~\label{minex}$\nu$ is not minimal since $\nu\geq 0$ with
$P(\nu>0)>0$ and $U_0(\nu)=0$.\er

\subsection{A solution to  BSDE ($g$, $\nu$)}
In the following, we find a bounded solution of BSDE ($g$, $\nu$)
where $0<\nu\leq T$ is a stopping time. Of course we can then
construct infinitely many bounded solutions for the BSDE.

{\bf Step 1.} \emph{For any $y_0<0$, construct an $\f-$predictable
process $H$ which can be dominated by $t\wedge \nu(\o)$ and
$t\wedge\nu(\o) +(1-\frac{t}{T})y_0$ for any $t$ small enough.}

Since $g$ is superquadratic and continuous,  we can get an
increasing sequence $\{x_i\}_{i=1}^{\infty}$ such that
$$ g(x_i)= i^{2}x_i^2$$
for any $i\geq \sqrt{a}+1$ where $a=\inf_{|x|>0}\frac{g(x)}{x^2}$
and $x_i=1$ for any $i< \sqrt{a}+1$.

Set $k_t=x_i$ when $t\in (\sum_{n=i+1}^{\infty}\delta_n,
\sum_{n=i}^{\infty}\delta_n]$ for $i\geq 1$ where
 $\delta_n=\frac{a_1}{g(x_n)n^2}$, and
$a_1=\frac{T}{\sum_{n=1}^{\infty}\frac{1}{g(x_n)n^2}}$ is such that
$\sum_{n=1}^{\infty}\d_n=T$.

We then have
$$\int_0^T(k_u)^2\,du=\sum_{i=1}^{\infty}x_i^2\frac{a_1}{g(x_i)i^2}<\infty,$$
$$\int_0^Tg(k_u)\,du=\sum_{i=1}^{\infty}g(x_i)\frac{a_1}{g(x_i)i^2}<\infty.$$
 Put $H_t=y_0+\int_0^tg(k_u)\,du-\int_0^tk_udB_u$.
\bl There exists $\Omega^{*}\subseteq \Omega$ with $P(\Omega^{*})=1$
satisfying for any $\varepsilon>0$ and $\o\in \Omega^{*}$ there is
$t_{\varepsilon}(\o)$ such that, for any $t<t_{\varepsilon}(\o)$,
$t\wedge \nu(\o)+(1-\frac{t}{T})y_0<H_t(\o).$\el

\emph{Proof. } It follows from the law of the iterated logarithm of
Brownian motion that there exists a set $\Omega^{*}$,
$P(\Omega^{*})=1$ satisfying:
 for any $\varepsilon>0$ and $\o\in \Omega^{*}$ there is
$t_{\varepsilon}(\o)<\nu$ such that $$\[\int_0^tk_udB_u\](\o)\leq
(1+\varepsilon)\sqrt{2\left(\int_0^tk_u^2\,du\right)\log\log\left(1/\int_0^tk_u^2\,du\right)}$$
for any $t<t_{\varepsilon}(\o)$.

Set
$F(t)=\int_0^tg(k_u)\,du-(1+\varepsilon)\sqrt{2\left(\int_0^tk_u^2\,du\right)\log\log\left(1/\int_0^tk_u^2\,du\right)}-\left(\frac{T-y_0}{T}\right)t$.
Now we want to prove $F(t)>0$ for $t$ small enough. Calculating the
differential of $F$ with respect to $t$, we have, for sufficiently
small $t$,
$$F'(t)> g(k_t)-1+\frac{y_0}{T}-\gamma k_t^2 \Big(\frac{1}{c_t}\log\log\frac{1}{c_t}\Big)^{1/2} $$
where $\gamma=(1+\varepsilon)2^{-\frac{1}{2}}$ and
$c_t=\int_0^tk_u^2\,du$. For $t\in (\sum_{n=i+1}^{\infty}\delta_n,
\sum_{n=i}^{\infty}\delta_n]$, we get
$$F'(t)> g(x_i)-1+\frac{y_0}{T}-\gamma x_i^2 \Big(\frac{1}{c_t}\log\log\frac{1}{c_t}\Big)^{1/2} $$
and \bea*c_t&\geq& \int_0^{\sum_{n= i+1}^{\infty}\delta_n}k_u^2\,du\\
&=&\sum_{n= i+1}^{\infty}x_n^2\delta_n\\
&=&a_1\sum_{n= i+1}^{\infty}n^{-4},\eea* for $i$ big enough. Thus
$$F'(t)> i^2x_i^2-1+\frac{y_0}{T}-\gamma x_i^2 \(\frac{1}{a_1\sum_{n= i+1}^{\infty}n^{-4}}\log\log\frac{1}{a_1\sum_{n=
i+1}^{\infty}n^{-4}}\)^{1/2}.$$ It follows from
$${ \lim_{i \rightarrow \infty}}\frac{i^2}{\gamma \(\frac{1}{a_1\sum_{n= i+1}^{\infty}n^{-4}}\log\log\frac{1}{a_1\sum_{n=
i+1}^{\infty}n^{-4}}\)^{1/2}}=+\infty$$ that there exist $0<t_0<T$
such that for any $t<t_0$, $F'(t)>0$. Since $F(0)=0$, we have, for
any $t<t_0$, \be~\label{fgeq0} F(t)>0.\ee

Thus for any $\o\in \Omega^{*}$ and $0<t<t_0\wedge
t_{\varepsilon}(\o)$, we have $$H_t(\o)-\Big[t\wedge
\nu(\o)+\big(1-\frac{t}{T}\big)y_0\Big]\geq F(t)>0.$$. \endpf

{\bf Step 2.} Since $H_t$ and $t\wedge \nu(\o)$ are $\{\f_t\}_{t\geq
0}$-predictable, we can define stopping times:
$$\tau_1^1\triangleq\inf\Big\{t>0\,\Big|\,H_t\leq t\wedge \nu+\big(1-\frac{t}{T}\big)y_0\Big\}\wedge T,$$
$$\tau_1^2\triangleq\inf\{t>0\,|\,H_t\geq t\wedge\nu\}\wedge T.$$
Define a random time
\begin{eqnarray}
\tau_1\triangleq
1_{\{\tau_1^1<\tau_1^2\}}\tau_1^1&+&1_{\{\tau_1^1\geq\tau_1^2\}}1_{\{\tau_1^2\leq
\nu+(1-\frac{\nu}{T})y_0\}}\frac{\tau_1^2-y_0}{T-y_0}T \nonumber \\
&+&1_{\{\tau_1^1\geq\tau_1^2\}}1_{\{\nu>\tau_1^2>
\nu+(1-\frac{\nu}{T})y_0\}}\frac{\nu-\tau_1^2+y_0}{y_0}T \nonumber \\
&+&1_{\{\tau_1^1\geq\tau_1^2\}}1_{\{\tau_1^2\geq \nu\}}T.
\end{eqnarray}
It is easy to verify that for any $\o\in \{\tau_1^1\geq\tau_1^2\}$,
$\tau_1(\o)\wedge \nu+(1-\frac{\tau_1(\o)}{T})y_0=\tau_1^2\wedge
\nu$.

\bl~\label{taustop} $\tau_1$ is a stopping time.\el

 \emph{Proof.} This is straightforward but for completeness we give a proof. The random time $\tau_1$ is defined by four parts without any
intersections. For the first part, it is easily verified by
\be~\label{meas1}\{1_{\{\tau_1^1<\tau_1^2\}}\tau_1^1\leq
t\}=\{\tau_1^1<\tau_1^2\}\cap \{\tau_1^1\leq t\}\in \f_t.\ee For the
second part, it is necessary  to check that
$$\{\tau_1^1\geq\tau_1^2\}\cap \big\{\tau_1^2\leq \nu+\big(1-\frac{\nu}{T}\big)y_0\big\}\cap \big\{\frac{\tau_1^2-y_0}{T-y_0}T\leq t\big\}\in \f_t.$$

It follows from $t(1-\frac{y_0}{T})+y_0\leq t$ that
$$\big\{\frac{\tau_1^2-y_0}{T-y_0}T\leq t\big\}=\big\{\tau_1^2\leq t(1-\frac{y_0}{T})+y_0\big\}\in \f_t$$
which implies that $\frac{\tau_1^2-y_0}{T-y_0}T$ is a stopping time.

If $\o \in \{\frac{\tau_1^2-y_0}{T-y_0}T\leq t\}$, then we have
$$\tau_1^2(\o)\leq t+(1-\frac{t}{T})y_0\le t.$$

Thus \begin{eqnarray}&&\{\tau_1^1\geq\tau_1^2\}\cap
\big\{\tau_1^2\leq \nu+\big(1-\frac{\nu}{T}\big)y_0\big\}\cap
\big\{\frac{\tau_1^2-y_0}{T-y_0}T\leq t\big\}
\nonumber\\
&=& \{\tau_1^1\geq\tau_1^2\}\cap\{\tau_1^2\leq t\}\cap
\big\{\frac{\tau_1^2-y_0}{T-y_0}T\leq \nu\big\}\cap
\big\{\frac{\tau_1^2-y_0}{T-y_0}T\leq t\big\} \nonumber\\
&=& \big[\{\tau_1^1\geq\tau_1^2\}\cap\{\tau_1^2\leq t\}\big]\cap
\big\{\frac{\tau_1^2-y_0}{T-y_0}T\leq \nu\wedge t\big\} \nonumber\\
&\in& \f_t. \end{eqnarray}

For the third part, observe that for any $\o \in \{\tau_1^2>
\nu+(1-\frac{\nu}{T})y_0\}\cap \{\frac{\nu-\tau_1^2+y_0}{y_0}T\leq
t\}$, we have
$$\nu(\o)-\tau_1^2(\o)\geq \big(\frac{t}{T}-1\big)y_0\geq 0.$$
Combining with $\nu(\o)-\tau_1^2(\o)< (\frac{\nu(\o)}{T}-1)y_0$, we
get $$\tau_1^2(\o)\leq \nu(\o)< t.$$ Therefore,
\begin{eqnarray}~\label{meas2}&&\{\tau_1^1\geq\tau_1^2\}\cap
\big\{\nu>\tau_1^2>
\nu+\big(1-\frac{\nu}{T}\big)y_0\big\}\cap
\big\{\frac{\nu-\tau_1^2+y_0}{y_0}T\leq t\big\} \nonumber\\
&=&\{\tau_1^1\geq\tau_1^2\}\cap \big\{\tau_1^2>
\nu+\big(1-\frac{\nu}{T}\big)y_0\big\}\cap
\big\{\frac{\nu-\tau_1^2+y_0}{y_0}T\leq t\big\} \nonumber\\
&=& \(\{\tau_1^1\geq\tau_1^2\}\cap \{\tau_1^2< t\}\)\cap\(
\big\{\frac{\tau_1^2-y_0}{T-y_0}T>\nu\big\}\cap\{\nu<
t\}\)\nonumber\\&&\cap\(
\big\{\tau_1^2+\big(\frac{t}{T}-1\big)y_0\leq \nu\big\}\cap\{\nu< t\}\) \nonumber\\
&\in&\f_t. \end{eqnarray} The fourth part is obviously
$\f_t$-measurable. Thus from (\ref{meas1}) to (\ref{meas2}) we get
that $\tau_1$ is a stopping time. \endpf

Define the predictable process $Z$ on the set $\{t\leq \tau_1\}$ as:
\be Z_t1_{\{t\leq \tau_1\}}=k_t1_{\{t\leq \tau_1^1\wedge
\tau_1^2\}}.\ee

\bl ~\label{Xline}Set $X_t\triangleq
y_0+\int_0^tg(Z_u)\,du-\int_0^tZ_udB_u$. We have
$X_{\tau_1}=\tau_1\wedge \nu+(1-\frac{\tau_1}{T})y_0$.\el

\emph{Proof.} The definitions of the stopping times yield that \bea*
X_{\tau_1}&=&  y_0+\int_0^{\tau_1}g(Z_u)\,du-\int_0^{\tau_1}Z_udB_u \\
&=&
y_0+\int_0^Tg(Z_u1_{\{t\leq\tau_1\}})\,du-\int_0^TZ_u1_{\{t\leq\tau_1\}}\,dB_u\\
&=&
y_0+\int_0^Tg(k_u)1_{\{t\leq\tau_1^1\wedge \tau_1^2\}}\,du-\int_0^Tk_u1_{\{t\leq\tau_1^1\wedge\tau_1^2\}}\,dB_u\\
&=& \(y_0+\int_0^{\tau_1^1}g(k_u)\,du-\int_0^{\tau_1^1}k_udB_u\)1_{\{\tau_1^1<\tau_1^2\}}\\
&&+\(y_0+\int_0^{\tau_1^2}g(k_u)\,du-\int_0^{\tau_1^2}k_udB_u\)1_{\{\tau_1^1\geq\tau_1^2\}}\\
&=&
\Big(\tau_1^1\wedge\nu+\big(1-\frac{\tau_1^1}{T}\big)y_0\Big)1_{\{\tau_1^1<\tau_1^2\}}+(\tau_1^2\wedge\nu)
1_{\{\tau_1^1\geq\tau_1^2\}} \\ &=&
\Big(\tau_1\wedge\nu+\big(1-\frac{\tau_1}{T}\big)y_0\Big)1_{\{\tau_1^1<\tau_1^2\}}+\Big(\tau_1\wedge\nu+\big(1-\frac{\tau_1}{T}\big)y_0\Big)
1_{\{\tau_1^1\geq\tau_1^2\}} \\
&=&\tau_1\wedge\nu+\big(1-\frac{\tau_1}{T}\big)y_0 \eea* which
completes the proof.

 {\bf Step 3.} \emph{Consider the solution in the time interval $(\tau_1,T)$.}

Construct $H_t^2=\tau_1\wedge
\nu+(1-\frac{\tau_1}{T})y_0+\int_{\tau_1}^tg(k_{u-\tau_1})\,du-\int_{\tau_1}^tk_{u-\tau_1}\,dB_u$
for any $t>\tau_1$.

Set $t=\tau_1+s$ where $s>0$. We have \begin{eqnarray}~\label{pt}
&&H_t^2-\Big(t\wedge
\nu+\big(1-\frac{t}{T}\big)y_0\Big)\nonumber\\
&=&\tau_1\wedge \nu+\big(1-\frac{\tau_1}{T}\big)y_0-\Big(t\wedge
\nu+\big(1-\frac{t}{T}\big)y_0\Big)+\int_{\tau_1}^tg(k_{u-\tau_1})\,du-\int_{\tau_1}^tk_{u-\tau_1}\,dB_u\nonumber\\
&\geq&\int_{0}^sg(k_{u})\,du-\int_{\tau_1}^{\tau_1+s}k_{u-\tau_1}\,dB_u+\frac{s}{T}y_0-s.\end{eqnarray}
Applying the law of the iterated logarithm of Brownian motion to
(\ref{pt}), we get that there is a set $\Omega^{*}\in \Omega$ with
$P(\Omega^{*})=1$ such that for any $\varepsilon>0$ and $\o \in
\Omega^{*}$ there exists a $s_{\varepsilon}(\o)$ satisfying for all
$s<s_{\varepsilon}(\o)$,
\bea*&&\int_{0}^sg(k_{u})\,du-\int_{\tau_1}^{\tau_1+s}k_{u-\tau_1}\,dB_u+\frac{s}{T}y_0-s
\\
&\geq&\int_{0}^sg(k_{u})\,du-(1+\varepsilon)\sqrt{2\left(\int_{\tau_1}^{\tau_1+s}k_{u-\tau_1}^2\,du\right)
\log\log\left(1/\int_{\tau_1}^{\tau_1+s}k_{u-\tau_1}^2\,du\right)}+\frac{s}{T}y_0-s\\
&=&\int_{0}^sg(k_{u})\,du-(1+\varepsilon)\sqrt{2\left(\int_{0}^{s}k_{u}^2\,du\right)
\log\log\left(1/\int_{0}^{s}k_{u}^2\,du\right)}+\frac{s}{T}y_0-s\\
&=& F(s). \eea* It follows from (\ref{fgeq0}) and (\ref{pt}) that
\be~\label{p2geq}H_{\tau_1+s}^2(\o)>(\tau_1+s)\wedge\nu(\o)+\big(1-\frac{\tau_1+s}{T}\big)y_0\ee
for all $0<s<s_{\e}(\o)\wedge t_0$.

Then similarly we can define stopping times: \be\tau_2^1\triangleq
\inf\Big\{t>\tau_1\,\Big|\,H_t^2\leq t\wedge \nu
+\big(1-\frac{t}{T}\big)y_0\Big\}\wedge T,\ee
\be\tau_2^2\triangleq\inf\Big\{t>\tau_1\,\Big|\,H_t^2\geq
t\wedge\nu\Big\}\wedge T,\ee and a random time
\begin{eqnarray}
\tau_2\triangleq
1_{\{\tau_2^1<\tau_2^2\}}\tau_2^1&+&1_{\{\tau_2^2\leq\tau_2^1\}}1_{\{\tau_2^2\leq
\nu+(1-\frac{\nu}{T})y_0\}}\frac{\tau_2^2-y_0}{T-y_0}T\nonumber\\
&+&1_{\{\tau_2^2\leq\tau_2^1\}}1_{\{\nu>\tau_2^2>
\nu+(1-\frac{\nu}{T})y_0\}}\frac{\nu-\tau_2^2+y_0}{y_0}T\nonumber \\
&+&1_{\{\tau_2^2\leq\tau_2^1\}}1_{\{\tau_2^2\geq \nu\}}T,
\end{eqnarray}
which is also a stopping time by a similar proof of $\tau_1$ in
Lemma \ref{taustop}.

{\bf Step 3.} \emph{Define the random times by transfinite
induction.}

The random time $\tau_{\alpha}$ for some ordinal number $\alpha$ is
defined by the following rules:
\begin{enumerate}
    \item $\tau_0=0$;
    \item If $E[\tau_{\a}]<T$, define $\tau_{\a+1}^1\triangleq
\inf\{t>\tau_{\a}\mid H_t^{\a+1}\leq t\wedge \nu
+(1-\frac{t}{T})y_0\}\wedge T$,
$$\tau_{\a+1}^2\triangleq\inf\{t>\tau_{\a}\mid H_t^{\a+1}\geq
t\wedge\nu\}\wedge T \mbox{ and }$$
\begin{eqnarray}
\tau_{\a+1}\triangleq
1_{\{\tau_{\a+1}^1<\tau_{\a+1}^2\}}\tau_{\a+1}^1&+&
1_{\{\tau_{\a+1}^2\leq\tau_{\a+1}^1\}}1_{\{\tau_{\a+1}^2\leq
\nu+(1-\frac{\nu}{T})y_0\}}\frac{\tau_{\a+1}^2-y_0}{T-y_0}T\nonumber\\
&+&1_{\{\tau_{\a+1}^2\leq\tau_{\a+1}^1\}}1_{\{\nu>\tau_{\a+1}^2>
\nu+(1-\frac{\nu}{T})y_0\}}\frac{\nu-\tau_{\a+1}^2+y_0}{y_0}T\nonumber \\
&+&1_{\{\tau_{\a+1}^2\leq\tau_{\a+1}^1\}}1_{\{\tau_{\a+1}^2\geq
\nu\}}T
\end{eqnarray}
where $H_t^{\a+1}=\tau_{\a}\wedge
\nu+(1-\frac{\tau_{\a}}{T})y_0+\int_{\tau_{\a}}^tg(k_{u-\tau_{\a}})\,du-\int_{\tau_{\a}}^tk_{u-\tau_{\a}}\,dB_u$
for $t>\tau_{\a}$.

    \item If $\beta$ is a limit number and satisfies $E[\tau_{\a}]<T$,
for all $\alpha<\beta$, then $\tau_{\b}\triangleq
\lim_{\a<\b}\tau_{\a}$.
\end{enumerate}
We adopt the symbol $\omega_1$ for the first uncountable ordinal and
let {\bf $\mathcal{O}$} be the well ordered set of all countable
ordinals, i.e. ordinals $\alpha<\omega_1$. Define

\be~\label{deflambda}\Lambda \triangleq\{\a\in{\bf \mathcal{O}}\mid
E[\tau_{\xi}]<T, \mbox{ for all }\xi <\a\}.\ee

 Since $\{E[\tau_{\a}]\}_{\a\in \Lambda}$ is
strictly increasing, $\Lambda$ is countable and hence there must
exist $\beta_0$ with $E[\tau_{\beta_0}]=T$, hence
$\tau_{\beta_0}=T$.

Define the predictable process $Z$ by \be Z_t=\sum_{0\leq
k<\b_0}k_{t-\tau_{k}}1_{\{\tau_{k}<t\leq \tau_{k}+\tau_{k+1}^1\wedge
\tau_{k+1}^2\}}\ee  and the stochastic process \be
X_t=y_0+\int_0^tg(Z_u)\,du-\int_0^tZ_udB_u \ee for any $t\in[0,T]$.

Similarly as lemma (\ref{Xline}), we have
$$X_{\tau_k}=\tau_k\wedge\nu+(1-\frac{\tau_k}{T})y_0$$
for any $k<\b_0$. Letting $k$ tend to $\b_0$, we get $X_T=\nu$.

Therefore, we constructed a solution $(X,Z)$ for the BSDE ($g$,$\nu$)
with $X_0=y_0<0$.

\begin{example}~\label{mincont} In this example, our goal is to construct
a bounded random variable $\zeta$ such that $U(\zeta)$ is a solution
of BSDE ($g$,\,$\zeta$) and $\zeta$ is not minimal when $g$ is
superquadratic. Define $\zeta=H_{\tau_1^1\wedge \tau_1^2}$, then
$(H^{\tau_1^1\wedge \tau_1^2}_t,k_t1_{\{t\leq \tau_1^1\wedge
\tau_1^2\}})_{0\leq t\leq T}$ is a solution to the
BSDE($g$,\,$\zeta$). It follows that
$$U(\zeta)=H^{\tau_1^1\wedge
\tau_1^2}.$$ Indeed, for any $t \in (0,T]$, set the probability
measure $Q^t$ via $\frac{dQ^t}{dP}=\mathcal{E}(q^t\cdot B)_T$ with
$q_s^t=g'(k_s)1_{\{t<s<\tau_1^1\wedge \tau_1^2\}}$. We then have
$$E\left[\exp\left(\frac{1}{2}\int_0^T|q^t_s|^2\,ds\right)\right]\leq\exp\left(\frac{1}{2}\int_t^T|g'(k_s)|^2\,ds\right)<\infty,$$
which implies $Q^t \sim P$. We deduce that
$$
U_t(\zeta)=E_{Q^t}\left[\zeta+\int_t^Tf(q^t_s)\,ds\bigg|\f_t\right]=H_t^{\tau_1^1\wedge
\tau_1^2}, \mbox{ for any } t\in(0,T],
$$
by the same argument of (\ref{ueqy}). Since
$H_{\cdot}^{\tau_1^1\wedge \tau_1^2}$ is continuous and
$U_{\cdot}(\zeta)$ is c\`{a}dl\`{a}g, we get $U_0(\zeta)=y_0$.

However, $\zeta$ is not minimal since $\zeta\geq y_0$ with
$P(\zeta>y_0)>0$ and $U_0(y_0)=y_0$.\endpf
\end{example}

\section{Existence of solution to BSDEs in the Markovian case}
From the last section, we know that the BSDE with superquadratic growth
is ill-posed. However we will show that in some particular Markovian case,
there exists a solution for such a BSDE.
\\ Define the diffusion
process $X^{t,x}$ to be the solution to the following  SDE:
\begin{eqnarray*}dX_s &=& b(s,X_s)\,ds+\sigma\,dB_s,\ \ t\leq s\leq
T,\\
X_s &=& x,\ \ 0\leq s\leq t,\end{eqnarray*} where $b:[0,T]\times \mathbb{R}^n\rightarrow \mathbb{R}^n$ is continuously differentiable with respect  to $x$ with bounded derivative $b_x$, and $\sigma
:[0,T]\rightarrow \mathbb{R}^{n \times d}$ is a constant (matrix).

Let us consider  BSDE (\ref{BSDEgsq}) with $\xi=\F (X_T^{t,x})$:

\be~\label{bsde_markov} Y_s=\F
(X_T^{t,x})-\int_s^Tg(Z_r)\,dr+\int_s^T Z_r\,dB_r,\quad s\in [0,T],
\ee where $g: R^{ d}\rightarrow R_{+}$ is a continuously
differentiable convex function with $g(0)=0.$ We suppose it is
superquadratic $\overline{\lim}_{|z|\rightarrow
\infty}\frac{g(z)}{|z|^2}=
    \infty.$ $f:R^{
d}\rightarrow R_{+}\cup\{\infty\}$ is the Fenchel-Legendre transform
of $g$: $$f(x)=\sup_{z\in R^d}(zx-g(z)),$$ then $f$ is also convex
and $f(0)=0.$


\subsection{Lipschitz case}
Let us first  consider the case when $\F$ is sufficiently smooth.
 \bt~\label{exist-unique} Suppose that $\Phi$ is bounded and Lipschitz. Then there exists a unique solution $(Y^{t,x},Z^{t,x})$ to
 BSDE (\ref{bsde_markov}) such that both processes $Y^{t,x}$ and $Z^{t,x}$
are bounded. Furthermore, the solution is a dynamic utility function
of the following form\begin{eqnarray}\label{sol-duf}
& &
Y_s^{t,x}\\
&=&\inf\bigg\{E_Q\Big[\F(X_T^{t,x})+\int_s^Tf(q_u)\,du\Big|\f_s\Big]\,\bigg|\,Q\sim P \bigg\}\nonumber\\
&=&\inf\bigg\{E_Q\Big[\F(X_T^{t,x})+\int_s^Tf(q_u)\,du\Big|\f_s\Big]\,\bigg|\,Q\sim
P, \nonumber \\&&\ \ \ \ \ E_Q\Big[\int_r^Tf(q_u)\,du\Big|\f_r\Big]
\leq 2\parallel\F\parallel_{\infty},\forall
r\in[0,T]\bigg\}\nonumber
\end{eqnarray} for any $s\in[0,T]$.\et
\emph{Proof.} First, let us suppose that $\F\in \mathcal{C}^1$ and
that
 $\F_x$ is bounded. We apply a
truncation argument to prove the existence of solution. Let us introduce the truncation function: for an
integer $N$, $\rho_N: \mathbb{R}^{1\times d} \rightarrow \mathbb{R}^+$
is smooth, such that $\forall |z|\leq N$, $\rho_N(z)=1$; and
$\forall |z|\geq N+1$, $\rho_N(z)=0$. Then it is obvious that $\rho_Ng$
is a bounded Lipschitz function. Hence for any $N$, there exists a unique
solution $(Y^{N;t,x},Z^{N;t,x})$ to the following BSDE:
 \be~\label{bsde_markov_N} Y_s=\F(X_T^{t,x})-\int_s^T(\rho_Ng)(Z_r)dr+\int_s^T Z_r\,dB_r.\ee
On the other hand, we  denote by $(F^{N;t,x},V^{N;t,x})$ the unique
solution to the following BSDE:

 \be~\label{U^Nequation}
F_s=\F_x(X_T^{t,x})\nabla_x
X_T^{t,x}-\int_s^T(\rho_Ng)_z(Z^{N;t,x}_r)V_r\,dr+\int_s^T
V_r\,dB_r,\ee where $\int_s^T V_r\,dB_r$ means$$\sum_{1\leq i\leq
d}\int_s^T V_r^{i}\,dB_r^i,$$ with $V^{i}$ denoting the $i$-th line
of the $d\times n$ matrix process $V$.

 We then have (see, e.g., \cite{PP92}):
\be
~\label{Z_U}Z^{N;t,x}_s=-F_s^{N;t,x}(\nabla_xX_s^{t,x})^{-1}\sigma.\ee
As for any $N$, $(\rho_Ng)_z(Z^{N;t,x})$ is bounded, we can apply a
Girsanov transformation to get:

 \be
F_s^{N;t,x}=\F_x(X_T^{t,x})\nabla_x X_T^{t,x}+\int_s^T
V^{N;t,x}_r\,dB^{N;t,x}_r,\ee where $B^{N;t,x}$ is a Brownian Motion
under an equivalent probability measure $Q^{N;t,x}$. Taking the
conditional expectation with respect to the measure $Q^{N;t,x}$, one
finally deduces  that $$|F_s^{N;t,x}|\leq
\parallel\F_x\parallel_{ \infty}\cdot\parallel\nabla_xX_T^{t,x}\parallel_{\infty}
$$
which implies that
\begin{eqnarray}|Z_s^{N;t,x}|&=&|F_s^{N;t,x}(\nabla_xX_s^{t,x})^{-1}\sigma|\nonumber \\
&\leq&\parallel\sigma\parallel\cdot
\parallel\F_x\parallel_{
\infty}\cdot e^{2\parallel b_x\parallel_{\infty}T}\\
\nonumber &:=& c.\nonumber\end{eqnarray}

The same argument (recall that $g(0)=0$) gives us also
that $$|Y_s^{N;t,x}|\leq \parallel\F\parallel_{\infty}.$$ Taking $N
\geq c$, then the solution $(Y^{N;t,x},Z^{N;t,x})$ to BSDE
(\ref{bsde_markov_N}) is actually a solution to BSDE
(\ref{bsde_markov}).

In the case when $\F$ is bounded and Lipschitz, we can also prove,
by a standard approximation, that there exists a  bounded
solution $(Y^{t,x},Z^{t,x})$ with $|Z^{t,x}|\leq c$
with $c=\parallel\sigma\parallel\cdot L_{\F}\cdot e^{2\parallel b_x\parallel_{\infty}T}$ where $L_{\F}$ is the Lipschitz constant of $\F$.

It is routine to prove the uniqueness of the bounded solution
($Y^{t,x}$, $Z^{t,x}$) where $Z^{t,x}$ is also bounded.


Finally, as $g_z(Z^{t,x})$ is bounded, and
$$E_{Q^{N;t,x}}\Big[\int_r^Tf(g_z(Z_u^{t,x}))\,du\Big|\f_r\Big] \leq
2\parallel\F\parallel_{\infty},\mbox{ for }N\ge c,\ \forall
r\in[0,T].$$ We conclude that $Y^{t,x}$ is a dynamic utility
function of the form (\ref{sol-duf}).
 \endpf

 \br A new solution $(Y,Z)$ can be constructed by the same
technique as before  with the process $Z$ {\bf unbounded}. \er
 We
define \be\label{utx} u(t,x):=Y_t^{t,x},\ee where
$(Y^{t,x},Z^{t,x})$ is the unique bounded solution to
(\ref{bsde_markov}) with $Z^{t,x}$ bounded. Since $\Phi$ is
Lipschitz,  $(Y^{t,x},Z^{t,x})$ is also the unique bounded solution
to (\ref{bsde_markov_N}) with $N\geq c$. An important property is
that $u(t,x)$ is deterministic. \br It follows from the classical
result of Markovian BSDEs that
$$Y_s^{t,x}=Y_t^{t,x},\  Z_s^{t,x}=0,\  \mbox{ for } s<t.$$
Besides, $(Y^{t,x},Z^{t,x})$ has the Markov property:
$$Y_s^{t,x}=u(s,X_s^{t,x}), \ \mbox{ for }s\ge t.$$\er

Furthermore, we have a uniqueness, a stability theorem and a strict
comparison theorem for the BSDEs. Thus we get the following
proposition. \bp~\label{pde-lip-prop} Suppose that $\Phi$ is bounded
and Lipschitz,  then $u(t,x)$ defined by (\ref{utx}) is bounded and
continuous on $[0,T]\times\mathbb{R}^n$ and a viscosity solution to
the PDE: \be~\label{pde-lip} \left\{\ba{l}
u_t(t,x)+\frac{1}{2}\mbox{trace}\big(\sigma\sigma^Tu_{xx}(
t,x)\big)+u_x(t,x)b(t,x) -g(-u_x(t,x)\sigma)=0,\\u(T,x)=\Phi(x).
    \ea\right.
\ee  \ep

\subsection{A Priori estimates of $Z$}

Now we suppose that both $\F$ and $\F_x$ are bounded. Let us first
suppose that $b\equiv0$ and $n=d$, $\sigma$ is the identity to
explain our main idea. In this case, $X_T^{t,x}=x+B_T-B_t$. Then
equation (\ref{Z_U}) turns out to be
$$Z^{N;t,x}_s=-F^{N;t,x}_s.$$
On the other hand, BSDE (\ref{U^Nequation}) becomes:
\be \left\{\ba{l}dZ^{N;t,x}_s=-(\rho_Ng)_z(Z^{N;t,x}_s)V_s^{N;t,x}\,ds+V_s^{N;t,x}\,dB_s,\q\q 0\leq s\leq T;\\
 Z^{N;t,x}_T=-\F_x(x+B_T-B_t)
    .\ea\right.
\ee This gives the following framework (taking $N\ge \parallel\F_x\parallel_\infty$):

\be~\label{general-framework} \left\{\ba{l}
dY_s=g(Z_s)\,ds-Z_s\,dB_s;\\dZ_s=-g_z(Z_s)V_s\,ds+V_s\,dB_s;\\
 Y_T=\xi\in L^{\infty}(\f_T),\q \q Y\mbox{ bounded}
    ,\ea\right.
\ee where $E\left[\int_0^T|Z_r|^2\,dr\right]< +\infty$ and
$\int_0^T|V_r|^2dr< +\infty\q P$ a.s. Thus we get special second
order backward stochastic differential equations (see  \cite{CSTV}
for a definition).\bt In the framework (\ref{general-framework}),
suppose there is a solution and

    1) The probability measure $Q$ with $\frac{dQ}{dP}=\mathcal{E}(g_z(Z)
B)_T$ is  equivalent to $P$;

2) $Z$ is a $Q$-martingale.

We then have \be|Z_s|\leq
2\parallel\xi\parallel_{\infty}(T-s)^{-\frac{1}{2}}, \ \ \forall
s\in [0,T).\ee

 Furthermore, if $f(g_z(\cdot)): \mathbb{R}^{ d}\rightarrow
\mathbb{R}^+ $ is convex, we also have:  \be f(g_z(Z_s))\leq
2\parallel\xi\parallel_{\infty}(T-s)^{-1}, \ \ \forall s\in
[0,T).\ee \et

\emph{Proof.} Under the measure $Q$, we get \be~\label{bsde-fg'}
dY_s=-f(g_z(Z_s))\,ds-Z_s\,dB^Q_s,\ee where
$B^Q_s=B_s-\int_0^sg_z^T(Z_r)dr$ is a $Q$-Brownian Motion.

Since $Y$ is bounded and $\int_0^sf(g_z(Z_r))dr$ is an increasing
process, it follows from Lemma \ref{bmo}  that $\int_0^sZ_rdB^Q_r$ is a BMO martingale under the measure
$Q$:

$$E_Q\bigg[\int_s^T|Z_r|^2dr\bigg|\f_s\bigg]\leq 4\parallel\xi\parallel^2_{\infty},$$
which implies, by Jensen's inequality,
$$|Z_s|^2(T-s)\leq 4\parallel\xi\parallel^2_{\infty},$$
i.e.
$$|Z_s|\leq 2\parallel\xi\parallel_{\infty}(T-s)^{-\frac{1}{2}}.$$
It follows from equation (\ref{bsde-fg'}) that
\begin{eqnarray*}
&& E_Q\bigg[\int_s^Tf(g_z(Z_r))dr\bigg|\f_s\bigg] \\ &=& -E_Q[\xi-Y_s|\f_s] \\
&\leq& 2\parallel\xi\parallel_{\infty}.
\end{eqnarray*}

 If $f(g_z(\cdot))$
is convex, then applying Jensen's inequality, we get
$$f(g_z(Z_s))\leq
2\parallel\xi\parallel_{\infty}(T-s)^{-1}.$$ \endpf
 In fact the condition 1) in the theorem is a constraint to make the process $\{Z_t\}_{0\leq t\leq T}$
not grow so fast as we constructed in the non-uniqueness theorem. In
this case, the solution is unique. We have the following remark. \br
Suppose there is a bounded solution $Y$ and the probability measure
$Q$ with $\frac{dQ}{dP}=\mathcal{E}(g_z(Z) B)_T$ is  equivalent to
$P$, then the solution is unique and
$$Y_s=E_Q\bigg[\xi+\int_s^Tf(g_z(Z_r))dr\bigg|\f_s\bigg],\ \ 0\leq s\leq T.$$
 \er

Let us consider the original BSDE (\ref{bsde_markov}) again. Taking
$N\ge c$ as in the proof of Theorem \ref{exist-unique}, we deduce the
following ``general''
 framework:

\be \left\{\ba{l}
dY_s=g(Z_s)\,ds-Z_s\,dB_s,\q\q Y_T=\F(X_T);\\dF_s=g_z(Z_s)V_sds-V_s\,dB_s, \q F_T=\F_x(X_T)\nabla_xX_T;\\
Z_s=-F_s(\nabla_x X_s)^{-1}\sigma,
    \ea\right.
\ee where $\F$ and $\F_x$ are bounded. Under the probability measure $Q$,  $B^Q_s=B_s-\int^s_0g_z(Z_r)dr$ is a Brownian
Motion, and the $``$general" framework
becomes:
 \be~\label{general-framework with
q} \left\{\ba{l}
dY_s=-f(g_z(Z_s))\,ds-Z_s\,dB^Q_s,\q Y_T=\F(X_T);\\dF_s=-V_s\,dB^Q_s,\q\q\q\q\q\q\q F_T=\F_x(X_T)\nabla_xX_T;\\
 Z_s=-F_s(\nabla_x X_s)^{-1}\sigma.
    \ea\right.
\ee Recall that $$d\nabla_x X_s=b_x(X_s)\nabla_x X_sds,$$ from which
we deduce
$$d(\nabla_x X_s)^{-1}=-(\nabla_x X_s)^{-1}(d\nabla_x X_s)(\nabla_x X_s)^{-1}=-(\nabla_x X_s)^{-1}b_x(X_s)\,ds.$$
Applying It\^o's formula, we deduce
$$dZ_s=-(dF_s)(\nabla_x X_s)^{-1}\sigma+F_s(\nabla_x X_s)^{-1}b_x(X_s)\sigma \,ds.$$
We suppose that there exists a constant $\lambda\geq 0$ such that
\be~\label{assumption} \forall \eta\in \mathbb{R}^n,
|\eta^T\sigma\sigma^Tb_x^T(x)\eta|\leq \lambda|\eta^T\sigma|^2.\ee
We then have
\begin{eqnarray*}
d(\exp(\l s)Z_s)&=& \l \exp(\l s)Z_sds+\exp(\l s)F_s(\nabla_x
X_s)^{-1}b_x(X_s)\sigma \,ds+dM_s \\
&=&F_s^*(\l I-b_x(X_s))\sigma \,ds+dM_s,
\end{eqnarray*}
where $M$ is a $Q$-martingale and
$$F_s^*=-\exp(\l s)F_s(\nabla X_s)^{-1}.$$
Finally,
$$d|\exp(\l s)Z_s|^2=d\langle M\rangle_s+2[\l |F^*_s\sigma|^2-F_s^*\sigma \sigma^Tb_x^T(X_s)(F_s^*)^T]\,ds+dM_s^*,$$
where $M^*$ is a $Q$-martingale, hence $|\exp(\l s)Z_s|^2$ is a
$Q$-submartingale.

\bp ~\label{Z-inequality} Let us suppose that $\F$ is bounded and
Lipschitz, $b$ and $\sigma$ satisfy the assumption
(\ref{assumption}), and ($Y^{t,x}$, $Z^{t,x}$) is the unique bounded
solution to BSDE (\ref{bsde_markov}). Then there exists a constant
$c_1>0$ such that \be |Z_s^{t,x}|\leq
c_1\parallel\F\parallel_{\infty}(T-s)^{-\frac{1}{2}},\ \  \forall s\in [0,T).\ee \ 
\ep

\emph{Proof.} First let us consider the smooth case when $\F$ and
$\F_x$ are bounded. Since $Z^{t,x}$ is of the framework
(\ref{general-framework with q}),
from Lemma \ref{bmo}, we have
$$E_Q\bigg[\int_s^T|Z_r^{t,x}|^2\,dr\bigg|\f_s\bigg]\leq 4 \parallel\F\parallel_{\infty}^2,$$
from which we deduce that
$$E_Q\bigg[\int_s^T\exp(2\l
r)|Z_r^{t,x}|^2\,dr\bigg|\f_s\bigg]\leq 4 \exp(2\l
T)\parallel\F\parallel_{\infty}^2.$$ As $\exp(2\l s)|Z_s^{t,x}|^2$
is a $Q$-submartingale, it follows that
$$\exp(2\l
s)|Z_s^{t,x}|^2(T-s)\leq 4 \exp(2\l
T)\parallel\F\parallel_{\infty}^2,$$ i.e.
$$|Z_s^{t,x}|\leq
c_1\parallel\F\parallel_{\infty}(T-s)^{-\frac{1}{2}},$$ where
$c_1=2\exp(\l T)$.

We can get the same estimate by a standard approximation when $\F$
is only bounded and  Lipschitz.\endpf

\br 
As an example, let us take $g(z)=|z|^q$ for $q\geq
 2$. Since $|\exp(\l s)Z_s^{t,x}|^2$ is a
$Q$-submartingale, it is clear that $|\exp(\l s)Z_s^{t,x}|^q$ is
also a $Q$-submartingale for $q\ge 2$. It follows from
(\ref{bsde-fg'}) that
 $$|Z_s^{t,x}|^q\leq C_q\parallel\F\parallel_\infty(T-s)^{-1}, s\in [t,T), $$
where $C_q>0$ is a constant depending only on $q$ and $\l$.
 Suppose $u$ is the bounded classical
solution to the following PDE: \be \left\{\ba{l}
u_t(t,x)+\frac{1}{2}\mbox{trace}\big(\sigma\sigma^T u_{xx}(
t,x)\big)+u_x(t,x)b(t) -|u_x(t,x)\sigma|^q=0,\\u(T,x)=\Phi(x).
    \ea\right.
\ee Since
$$Z_s^{t,x}=-u_x(s,X_s^{t,x})(\nabla_xX_s^{t,x})^{-1}\sigma,\
\forall s \in [t,T],$$ we deduce that
$$|u_x(t,x)\sigma|\leq (C_q\parallel\F\parallel_{\infty})^{1/q}(T-t)^{-1/q}.$$
The same type of estimate is given by Gilding et al. in \cite{ggk}
using Bernstein's technique, in the case when $b=0$ and $\sigma$ is
the identity.
 \er
\subsection{Lower semi-continuous case }

Notice that $Z^{t,x}$ is bounded when $\F$ is bounded and Lipschitz.
 The bound, however, depends on the Lipschitz constant. The advantage of the estimate in  Proposition \ref{Z-inequality} is that the estimate only depends
on $\parallel\F\parallel_{\infty}$. This allows us to weaken the
hypothesis further. \bp ~\label{lower-semi-case}Let us suppose that
$\F$ is bounded and lower semi-continuous, and $b$ and $\sigma$
satisfy the assumption (\ref{assumption}). Then there exists a
bounded solution ($\underline{Y}^{t,x}$, $\underline{Z}^{t,x}$) to
BSDE (\ref{bsde_markov}) such that
\begin{equation}\label{Z}|\underline{Z}_s^{t,x}|\le c_1\parallel\F\parallel_{\infty}(T-s)^{-\frac{1}{2}},\ \  \forall s\in [t,T).
\end{equation}
\ep

\emph{Proof.} For each integer $m\geq 0$, construct the function
$$\underline{\F}_m(u)=\inf\{\F(p)+m|p-u|: p\in \mathbb{R}^n\}.$$
Then $\underline{\F}_m$ is well defined and globally Lipschitz with Lipschitz
constant $m$. Moreover, $(\underline{\F}_m)_{m\geq 0}$ is increasing
and converges pointwise to $\F$ with
$$-\parallel\F\parallel_{\infty}\leq \underline{\F}_m\leq \F.$$
Let ($\underline{Y}^{m;t,x}$, $\underline{Z}^{m;t,x}$) be the
bounded solution to BSDE ($g$, $\underline{\F}_m(X_T^{t,x})$). It
follows from the classical comparison theorem that
$$-\parallel\F\parallel_{\infty}\leq \underline{Y}^{0;t,x}\leq \underline{Y}^{m;t,x}\leq \underline{Y}^{m+1;t,x}\leq \parallel\F\parallel_{\infty},$$
and from Proposition \ref{Z-inequality},
\begin{equation}\label{Zm}|\underline{Z}_s^{m;t,x}|\leq
c_1\parallel\F\parallel_{\infty}(T-s)^{-\frac{1}{2}},\quad s\in [0,T).
\end{equation}
 For any
fixed $T'\in (0,T)$, ($\underline{Y}^{m;t,x}$, $\underline{Z}^{m;t,x}$) satisfies
\begin{equation}\label{bsdemT'}\underline{Y}_s^{m;t,x}=\underline{Y}_{T'}^{m;t,x}-\int_s^{T'}(\rho_Mg)(\underline{Z}_r^{m;t,x})\,dr+\int_s^{T'}\underline{Z}_r^{m;t,x}\,dB_r,\q\q \forall s \in[0,T'],\end{equation}
where
$$M=c_1\parallel\F\parallel_{\infty}(T-T')^{-\frac{1}{2}}.$$
Moreover, by Lemma \ref{bmo},
$$E\left[\int_s^T |\underline{Z}_r^{m;t,x}|^2dr\Big|{\cal F}_s\right]\le 4\parallel\F\parallel_\infty^2.$$
The classical stability theorem (see N. El Karoui et al \cite{EPQ})
for Lipschitz generators implies
$$\lim_{m,m'\rightarrow\infty}E\left[\int_0^{T'}\left|\underline{Z}_r^{m;t,x}-\underline{Z}_r^{m';t,x}\right|^2dr\right]=0.$$
So define
$$\underline{Y}^{t,x}=\lim_{m\rightarrow\infty} \underline{Y}^{m;t,x},\q\q \underline{Z}^{t,x}=\lim_{m\rightarrow\infty} \underline{Z}^{m;t,x}.$$
Then by passing to the limit when $m\rightarrow \infty$ in (\ref{bsdemT'}), we conclude that
for any
fixed $T'\in (0,T)$, ($\underline{Y}^{t,x}$, $\underline{Z}^{t,x}$) satisfies
\begin{equation}\label{bsdeT'}{Y}_s={Y}_{T'}-\int_s^{T'}g({Z}_r)\,dr+\int_s^{T'}{Z}_r\,dB_r,\q\q \forall s \in[0,T'],
\end{equation}
and
$$|\underline{Y}_s^{t,x}|\le\parallel\F\parallel_\infty,\quad E\left[\int_s^T |\underline{Z}_r^{t,x}|^2dr\Big|{\cal F}_s\right]\le 4\parallel\F\parallel_\infty^2,\quad s\in [0,T].$$

On the other hand, we have $$\liminf_{s\rightarrow T}\underline{Y}_s^{t,x}\geq
\liminf_{s\rightarrow
T}\underline{Y}_s^{m;t,x}=\underline{\F}_m(X_T^{t,x}) \mbox{ for any
} m\in\mathbb{N}, \ P\ a.s.$$ which implies $\liminf_{s\rightarrow
T}\underline{Y}_s^{t,x}\geq{\F}(X_T^{t,x}), \ P\ a.s.$

Since $\underline{\F}_m$ is bounded and Lipschitz for any
$m\in\mathbb{N}$, it follows from Theorem \ref{exist-unique} that
$\underline{Y}_s^{m;t,x} \leq E[\underline{\F}_m(X_T^{t,x})|\f_s]$.
We then get
\begin{eqnarray*}
\limsup_{s\rightarrow
T}\underline{Y}_s^{t,x}&=&\limsup_{s\rightarrow
T}\lim_{m\rightarrow \infty}\underline{Y}^{m;t,x}_s \\
&\leq&\limsup_{s\rightarrow T}\lim_{m\rightarrow
\infty}E[\underline{\F}_m(X_T^{t,x})|\f_s] \\
&=& \limsup_{s\rightarrow T}E[{\F}(X_T^{t,x})|\f_s] \\
&=& {\F}(X_T^{t,x}),\ \  P\ a.s.
\end{eqnarray*}
Hence $\lim_{s\rightarrow T}\underline{Y}_s^{t,x}= {\F}(X_T^{t,x})$.

Finally, passing to the limit when $T'\rightarrow T$ in
(\ref{bsdeT'}), we conclude that ($\underline{Y}^{t,x}$,
$\underline{Z}^{t,x}$) is a bounded solution to BSDE
(\ref{bsde_markov}). By passing to the limit when
$m\rightarrow\infty$ in (\ref{Zm}), we derive (\ref{Z}) immediately.
\endpf

\subsection{Bounded and Continuous Case}

In the smooth case, the dynamic utility function is a solution to
BSDE (\ref{bsde_markov}) by Theorem \ref{exist-unique}. This remains
true in more general case.

\bp ~\label{contin-case}Let us suppose that $\F$ is bounded and
continuous, and $b$ and $\sigma$ satisfy the assumption
(\ref{assumption}). Then there exists a bounded solution
($\bar{Y}^{t,x}$, $\bar{Z}^{t,x}$) such that
$$\bar{Y}_s^{t,x}=\inf\Big\{E_Q\big[\F(X_T^{t,x})+\int_s^Tf(q_r)dr\big|\f_s\big]\,\Big|\,Q\sim P\Big\}.$$\ep

\emph{Proof.}  By the same technique as that used in Proposition
\ref{lower-semi-case}, let us define the function
$$\bar{\F}_m(u)=\sup\{\F(p)-m|p-u|: p\in \mathbb{R}^n\}$$
for each integer $m\geq 0.$ Then $\bar{\F}_m$ is also bounded and
globally Lipschitz  with Lipschitz constant $m$.
$(\bar{\F}_m)_{m\geq 0}$ is decreasing and converges pointwise to
$\F$ with
$$\parallel\F\parallel_{\infty}\geq \bar{\F}_m\geq \F.$$
Let ($\bar{Y}^{m;t,x}$, $\bar{Z}^{m;t,x}$) be the solution to
BSDE ($g$, $\bar{\F}_m(X_T^{t,x})$). It follows from the same
argument as that in Proposition \ref{lower-semi-case} that
by setting
$$\bar{Y}^{t,x}=\lim_{m\rightarrow
\infty}\bar{Y}^{m;t,x},\q\q\q\q\q\bar{Z}^{t,x}=\limsup_{m\rightarrow
\infty} \bar{Z}^{m;t,x},$$
$(\bar{Y}^{t,x},\bar{Z}^{t,x})$ satisfies (\ref{bsdeT'}). On the other hand, since
$$\underline{Y}_s^{m;t,x}\le \bar{Y}_s^{m;t,x}\le E[\bar{\F}_m(X_T^{t,x})|{\cal F}_s],$$
we deduce that
$$\underline{Y}_s^{t,x}\le \bar{Y}_s^{t,x}\le E[{\F}(X_T^{t,x})|{\cal F}_s],$$
which implies that $$\lim_{s\rightarrow T}\bar{Y}_s^{t,x}=\F(X_T^{t,x}).$$

Hence ($\bar{Y}^{t,x}$, $\bar{Z}^{t,x}$) is also a bounded solution to BSDE ($g$,
$\F(X_T^{t,x})$).  Lemma \ref{ugeqy} implies that any bounded
solution of BSDE is less than or equal to the corresponding dynamic utility
function,
$$\bar{Y}^{t,x}_s\leq \essinf\Big\{E_Q\big[\F(X_T^{t,x})+\int_s^Tf(q_r)dr\big|\f_s\big]\,\Big|\,Q\sim P\Big\}.$$
 Finally, as
$$\bar{Y}_s^{m;t,x}=\essinf\Big\{E_Q\big[\bar{\F}_m(X_T^{t,x})+\int_s^Tf(q_r)dr\big|\f_s\big]\,\Big|\,Q\sim P\Big\}$$
and $\bar{\F}_m(X_T^{t,x})$ converges decreasingly to
${\F}(X_T^{t,x})$, we deduce
$$\bar{Y}^{t,x}_s\geq
\essinf\Big\{E_Q\big[\F(X_T^{t,x})+\int_s^Tf(q_r)dr\big|\f_s\big]\,\Big|\,Q\sim
P\Big\}.$$ Combining the above, we conclude that the solution $\bar{Y}^{t,x}$ is a dynamic
utility function.\endpf

Notice that we used $\underline{\F}^m$ to approximate $\F$ in the
lower semi-continuous case. In the continuous case, we can show that
both $\underline{Y}^{m;t,x}$ and $\overline{Y}^{m;t,x}$  converge to
the same limit. Now first let us consider the uniformly continuous
case. \bt Suppose that $\F$ is bounded and uniformly continuous. We
then have
\be\bar{Y}^{t,x}_s=U_s(\F(X_T^{t,x}))=\underline{Y}^{t,x}_s.\ee\et

\emph{Proof.} It follows from the uniform continuity of $\F$ that
 both $(\bar{\F}_m)_{m\geq 0}$ and $(\underline{{\F}}_m)_{m\geq
0}$ converge to $\F$ with the norm
$\parallel\cdot\parallel_{\infty}$.

Indeed, the uniform continuity of $\F$ implies that, for any
$\epsilon >0$, there exists $\delta_\epsilon>0$, such that if
$|p-u|\leq \delta(\epsilon)$, then
$$|\F(p)-\F(u)|\leq \epsilon.$$
By the definition, we get
$$\F(u)-\underline{\F}_m(u)=\sup\{\F(u)-\F(p)-m|p-u|: p\in \mathbb{R}^n\}.$$
But we have
\begin{eqnarray*}
&&\F(u)-\F(p)-m|p-u| \\
&=&
(\F(u)-\F(p))1_{\{|p-u|<\delta_{\epsilon}\}}+(\F(u)-\F(p))1_{\{|p-u|\geq\delta_{\epsilon}\}}-m|p-u|
\\
&\leq& \epsilon
+2\frac{\parallel\F\parallel_{\infty}}{\delta_{\epsilon}}|p-u|-m|p-u|
\end{eqnarray*}
from which we deduce that if
$m>2\frac{\parallel\F\parallel_{\infty}}{\delta_{\epsilon}}$, then
$$0\leq \F(u)-\underline{\F}_m(u)\leq \epsilon,\forall u \in \mathbb{R}^n,$$
hence, $$\lim_{m\rightarrow
\infty}\parallel\F-\underline{\F}_m\parallel_{\infty}=0.$$

 Combining with the convergence of a dynamic utility function and
$$\underline{Y}_s^{m;t,x}=\essinf\Big\{E_Q\big[\underline{\F}_m(X_T^{t,x})+\int_s^Tf(q_r)dr\big|\f_s\big]\,\Big|\,Q\sim P\Big\},$$
we get
$$\underline{Y}^{t,x}_s=U_s(\F(X_T^{t,x}))=\essinf\Big\{E_Q\big[\F(X_T^{t,x})+\int_s^Tf(q_r)dr\big|\f_s\big]\,\Big|\,Q\sim P\Big\}.$$
By the same argument or simply by Proposition
\ref{contin-case}, we have
$$\bar{Y}^{t,x}_s=U_s(\F(X_T^{t,x}))=\underline{Y}^{t,x}_s.$$ \endpf

Let us now consider the general case: $\F$ is bounded and
continuous. \bt ~\label{same_limit}Suppose that $\F$ is bounded and
continuous and $f$ satisfies the assumption: there exists a constant
$M$ such that \be\alpha:=\min_{|x|=M}\{f(x)\}>0.\ee We then have
\be\bar{Y}^{t,x}_s=U_s(\F(X_T^{t,x}))=\underline{Y}^{t,x}_s,\ \
\forall s\in[0,T].\ee \et

\emph{Proof.} First, in the special case $X_T^{t,x}=x+B_T-B_t$, let
us consider $\overline{Y}^{m;t,x}-\underline{Y}^{m;t,x}$. It follows
from (\ref{sol-duf}) that
\begin{eqnarray*}
&&\overline{Y}_0^{m;t,x}-\underline{Y}_0^{m;t,x} \\
&=&\essinf\Big\{E_Q\big[\overline{\F}_m(X_T^{t,x})+\int_0^Tf(q_u)\,du\big]\,\Big|\,Q\sim
P, E_Q\big[\int_0^Tf(q_u)\,du\big] \leq
2\parallel\F\parallel_{\infty}\Big\}\\&&-\essinf\Big\{E_Q\big[\underline{\F}_m(X_T^{t,x})+\int_0^Tf(q_u)\,du\big]\,\Big|\,Q\sim
P,  E_Q\big[\int_0^Tf(q_u)\,du\big] \leq
2\parallel\F\parallel_{\infty}\Big\}\\
&\leq&
\esssup\Big\{E_Q\big[\overline{\F}_m(X_T^{t,x})-\underline{\F}_m(X_T^{t,x})\big]\,\Big|\,Q\sim
P, E_Q\big[\int_0^Tf(q_u)\,du\big] \leq
2\parallel\F\parallel_{\infty}\Big\}.
\end{eqnarray*}
Denoting $\overline{\F}_{m}-\underline{\F}_m$ as $\Psi_m$, then
$\Psi_m$ is continuous. We then analyze
\begin{eqnarray*}
\Psi_m(x+B_T-B_t) &=&\Psi_m\bigg(x+B^Q_T-B^Q_t+\int_t^Tq_udu\bigg),
\end{eqnarray*}
in three parts. First, we have, for $N>M$,
\begin{eqnarray*}&&E_Q\Big[\Psi_m\Big(x+B^Q_T-B^Q_t+\int_t^Tq_udu\Big)1_{\{
|\int_t^Tq_udu|>N\}}\Big]\\
&\leq&2\parallel\F\parallel_{\infty}Q\Big(\Big\{\big|\int_t^Tq_udu\big|>N\Big\}\Big)\\
&\leq&
2\parallel\F\parallel_{\infty}\frac{E_Q[|\int_t^Tq_udu|]}{N}.\end{eqnarray*}
Combining with \begin{eqnarray*} &&E_Q\Big[\Big|\int_t^Tq_udu\Big|\Big]\\
&\leq& E_Q\Big[\int_t^T\big|q_u\big|1_{\{|q_u|<
M\}}\,du\Big]+E_Q\Big[\int_t^T\big|q_u\big|1_{\{|q_u|\geq M\}}\,du\Big]\\
&\leq& MT+E_Q\Big[\int_0^T\frac{M}{\alpha}f(q_u)1_{\{|q_u|\geq M\}}\,\,du\Big]\\
&\leq&MT+2\frac{M}{\alpha}\parallel\F\parallel_{\infty},
\end{eqnarray*}
we deduce that
$$E_Q\Big[\Psi_m\Big(x+B^Q_T-B^Q_t+\int_t^Tq_udu\Big)1_{\{
|\int_t^Tq_udu|>N\}}\Big]\leq \frac{c_1}{N},$$ where $c_1>0$ is a
constant. Second, we have\begin{eqnarray*}
&&E_Q\Big[\Psi_m\Big(x+B^Q_T-B^Q_t+\int_t^Tq_udu\Big)1_{\{|B_T^Q-B^Q_t|>
N,
|\int_t^Tq_udu|\leq N\}}\Big]\\
&\leq&2\parallel\F\parallel_{\infty}\frac{E_Q[|B_T^Q-B^Q_t|]}{N}\\
&=&\frac{c_2}{N}
\end{eqnarray*}
where $c_2>0$ is a constant independent of $Q$ and $m$. Third, we
have
\begin{eqnarray*}
&&E_Q\Big[\Psi_m\Big(x+B^Q_T-B^Q_t+\int_t^Tq_udu\Big)1_{\{|B_T^Q-B^Q_t|\leq
N,
|\int_t^Tq_udu|\leq N\}}\Big]\\
&\leq&\sup_{|y|\leq |x|+2N}\Psi_m(y).
\end{eqnarray*}
It follows from the preceding three estimates that
$$E_Q[\Psi_m(x+B_T-B_t)]\leq \frac{c_1+c_2}{N}+\sup_{|y|\leq |x|+2N}\Psi_m(y),$$
which implies that
\begin{eqnarray*}
\limsup_{m\rightarrow
\infty}(\overline{Y}_0^{m;t,x}-\underline{Y}_0^{m;t,x})&\leq&\limsup_{m\rightarrow
\infty}\[\frac{c_1+c_2}{N}+\sup_{|y|\leq |x|+2N}\Psi_m(y)\]\\
&=&\frac{c_1+c_2}{N}.
\end{eqnarray*}
Since $\overline{Y}_0^{m;t,x}-\underline{Y}_0^{m;t,x}\geq 0$, by
letting $N$ tend to $\infty$, we deduce that
\be\overline{Y}_0^{t,x}=\underline{Y}_0^{t,x}.\ee Combining with the
Markov property of $\overline{Y}^{t,x}$ and $\underline{Y}^{t,x}$,
we conclude that $\overline{Y}_s^{t,x}=\underline{Y}_s^{t,x}$, for
any $0< s< T$.

Notice that essentially we have made use of the simple fact that
$$\lim_{m\rightarrow
\infty}\sup_{|y|\leq c}\Psi_m(y)=0,$$ where $c$ is a constant. So in
the general case,
$$X_s^{t,x}=x+\int_t^s b(u,X_u^{t,x})\,du+\sigma \left(B_s^Q-B^Q_t+\int_t^sq_u\,du\right),$$
since $b$ is Lipschitz, applying Gr{o}nwall's inequality, we get
$$|X_T^{t,x}|\leq\sup_{t \leq s\leq T} |X_s^{t,x}|\le C \left(1+ \int_t^T |q_u|\,du + \sup_{t \leq s\leq T} |B_s^Q-B^Q_t|\right),$$
where $C>0$ is a constant. Hence the same proof
works.\endpf

Now we define
 $$\overline{u}^m(t,x):=\overline{Y}^{m;t,x}_t,\ \underline{u}^m(t,x):=\underline{Y}^{m;t,x}_t,$$
and \be~\label{newutx} u(t,x):=U_t(\F(X_T^{t,x})).\ee We have the
following theorem.

\bt Suppose that $\F$ is bounded and continuous and $f$ satisfies
the assumption: there exists a constant $M>0$ such that
$$\alpha:=\min_{|x|=M}\{f(x)\}>0;$$ and that $b$ and $\sigma$ satisfy the assumption
(\ref{assumption}). Then $u(t,x)$ defined by (\ref{newutx}) is a
bounded and continuous deterministic function on
$[0,T]\times\mathbb{R}^n$ and it is a viscosity solution to PDE
(\ref{pde-lip}).\et

\emph{Proof.}  Theorem \ref{same_limit} implies that
$\{\overline{u}^m(t,x)\}_{m=1}^{\infty}$ (resp.
$\{\underline{u}^m(t,x)\}_{m=1}^{\infty}$) converges decreasingly
(resp. increasingly) to $u(t,x)$. Combining with the continuity of
 $\overline{u}^m(t,x)$ and
 $\underline{u}^m(t,x)$, we deduce that $u(t,x)$ is
 continuous.

 By Dini's theorem, they converge to $u$ uniformly
 in any compact set. This implies that $u$ is a viscosity solution by the stability theorem of viscosity solutions to
 PDEs
(see, e.g., \cite{CIL}) .
\endpf


\medskip

{\bf Acknowledgments} This research was
sponsored by a grant of Credit Suisse as well as by a grant
NCCR-Finrisk. The text only reflects the opinion of the authors.
We also thank M. Ben-Artzi  for stimulating discussions
regarding the connection with viscous Hamilton-Jacobi equations.

\end{document}